\documentclass[twocolumn]{svjour3}  
\smartqed  

\usepackage[dvipdfmx]{color}
\usepackage{graphicx}
\usepackage[T1]{fontenc}
\usepackage{mathptmx}  
\usepackage[scaled]{helvet}  
\usepackage{courier}
\usepackage[subrefformat=parens]{subcaption}

\makeatletter
\let\cl@chapter\undefined
\makeatletter
\usepackage{amsmath,amssymb}   
\usepackage{mathtools} 
\usepackage{cleveref}  
\Crefname{equation}{Eq.}{Eqs.}%
\Crefname{figure}{Fig.}{Figs.}%
\usepackage{amsbsy}
\usepackage{bm}
\usepackage{algorithmic}
\usepackage{algorithm}

\begin{document}

\title{Application of Linear Filter and Moment Equation for Parametric Rolling in Irregular Longitudinal Waves
}

\author{Yuuki Maruyama      \and Atsuo Maki \and
        Leo Dostal          \and Naoya Umeda
}

\institute{Yuuki Maruyama \and Atsuo Maki \and Naoya Umeda \at
              Osaka University, 2-1 Yamadaoka, Suita, Osaka, Japan \\
              \email{yuuki\_maruyama@naoe.eng.osaka-u.ac.jp} 
           \and
           Leo Dostal \at
           Institute of Mechanics and Ocean Engineering, Hamburg University of Technology, 21043 Hamburg, Germany
}

\date{Received: date / Accepted: date}

\maketitle

\begin{abstract}
Parametric rolling is one of the dangerous dynamic phenomena. In order to discuss the safety of a vessel when a dangerous phenomenon occurs, it is important to estimate the exceedance probability of certain dynamical behavior of the ship with respect to a certain threshold level. In this paper, the moment values are obtained by solving the moment equations. Since the stochastic differential equation(SDE) is needed to obtain the moment equations, the Autoregressive Moving Average(ARMA) filter is used. The effective wave is modeled by using the 6th-order ARMA filter. In addition, the parametric excitation process is modeled by using a linear-nonlinear transformation obtained from the relationship between GM and wave amplitude. The resulting system of equations is represented by the 8th-order Itô stochastic differential equation, which consists of a second-order SDE for the ship motion and a 6th-order SDE for the effective wave. This system has nonlinear components. Therefore, the cumulant neglect closure method is used as higher-order moments need to be truncated. Furthermore, the probability density function of roll angle is determined by using moment values obtained from the SDE and the moment equation. Here,  two types of the probability density function are suggested and have a good agreement.

\end{abstract}

\section{Introduction}
When roll motion in irregular waves is considered in the framework of the stochastic theory, the irregular excitation due to waves is assumed as a stochastic process. However, the real noise such as the processes of wind or sea waves are not of white noise type but of colored noise type. Therefore, a real noise process for sea waves should be generated from filtered white-noise via the stochastic differential equation~(SDE)~\cite{Dostal2011}. From such a point of view, Spanos~\cite{SPANOS1983}\cite{SPANOS1986}, Flower~\cite{Flower1983}~\cite{Flower1985}, and Thampi~\cite{THAMPI1992} applied earlier a linear filter for the generation of stochastic time series based on wave spectra (Pierson-Moskowitz wave spectrum and JONSWAP spectrum) using the ARMA process, which was a good approximation of the process obtained from sea wave spectra. On the other hand, concerning wind process generation, Nichita~\cite{Nichita2002}, Dostal et al.~\cite{Dostal2020_windwave}, and Maki et al.~\cite{maki2021wind} also utilized liner filters for white noise. 

In order to discuss the safety of a vessel when a dangerous phenomenon occurs, it is important to estimate the exceedance probability of certain dynamical behavior of the ship with respect to a certain threshold level. For this purpose, it is desirable to analytically or theoretically obtain a probability density function (PDF) of a critical dynamical response. One of the candidates to theoretically obtain the PDF in a nonlinear system is to solve the Fokker-Planck-Kolmogorov~(FPK) equation analytically~\cite{Caugy1963-1}. Unfortunately, the application of this method is slightly limited to  simple low dimensional nonlinear systems driven by white noise, such that the direct application of this method to the problem of parametric roll motion driven by colored noise is not possible. In order to approach such problems the equivalent linearization method was developed~\cite{Caugy1963-2}\cite{Socha2008}, which approximates a nonlinear system using linear components. 
One of the successive method to solve the FPK equation for this parametric rolling problem is the stochastic averaging method~\cite{Dostal2011}\cite{Roberts1982}\cite{Maruyama2021}. In this approach, the FPK equation for the averaged system is analytically solved. 
On the other hand, another method is to use moment equations~\cite{Bover1978}\cite{Wu1987}. In the case of using higher-order moment equations, a nonlinear effect can be reflected in the analysis.  As is well known, the response for a strong nonlinear system has a non-Gaussian distribution~\cite{YOKOYAMA2013}. In the case of the parametric rolling problem, the PDF of roll angle is not Gaussian as pointed out by Belenky~\cite{Belenky2011}.
So far, some researches use a combination of moment equations and a linear filter, which has been applied in the field of naval architecture and ocean engineering. For example, Francescutto et al.~\cite{FRANCESCUTTO2003} and Su et al.~\cite{Su2011} considered the roll motion in beam seas using a 4th-order linear filter and a moment equation. Chai et al.~\cite{Chai2016} analyzed the response of parametric rolling in irregular waves by using Monte Carlo simulation~(MCS). Dostal et al.~\cite{Dostal2011} used the Local Statistical Linearization in combination with moment equations.

In this study, the variation of metcentric height~(GM) is modeled by using the Autoregressive Moving Average~(ARMA) filter and the linear-nonlinear transformation. Then the results are examined by comparison with result obtained by using the superposition principle. The objective of this paper is to obtain the PDF of the roll angle. Thus, this paper is concerned with the moment equation derived from the stochastic differential equation representing the ship motion and the ocean wave elevation. Thereby, the PDF of the roll angle is determined by using stochastic moments. By this means, two types of PDF are chosen and their coefficients are optimized.
    
\section{Linear Filter}
    %
    When a stochastic differential equation~(SDE), is used  to model parametric rolling, then it is common to model the parametric excitation process as a stationary random process with zero mean~\cite{Dostal2011}\cite{Dostal2012}.
    It is generally acceptable that a parametric excitation (real noise) itself, calculated on the basis of hydrodynamic assumptions, is regarded as a stochastic process. 
    However, an It\^{o} stochastic differential equation, which uses a driving noise derived from white noise, is needed in order to determine moment equations. Therefore it is appropriate to represent the parametric excitation process approximately by a stochastic differential equation.
    In this study, as a means of approximating real noise from white noise, the time series of the effective wave is modeled using an Autoregressive Moving Average (ARMA) process. In addition, the parametric excitation process is determined by using a linear-nonlinear transformation obtained from the relationship between GM and wave amplitude as shown in Fig.\ref{fig:HBL_vs_GM}.

    In this study, the ITTC spectrum is used as the ocean wave spectrum, which is given by
     \begin{equation}
        S_{\mathrm{w}}(\omega) = \dfrac{173H_{1/3}^{2}}{T_{01}^{4}\omega^{5}} \exp\left(-\dfrac{691}{T_{01}^{4}\omega^{4}}\right)
    \end{equation}
    Here, $T_{01}$ and $H_{1/3}$ are the wave mean period and the significant wave height, respectively.
    In this study, the GM variation is calculated by applying the concept of Grim’s effective wave\cite{Grim1961} to irregular waves. The transfer function $H_{\zeta}^{}$ which is necessary to obtain the effective wave spectrum $S_{\mathrm{eff}}^{}$, can be represented as in the following equation\cite{Umeda1991}:
    \begin{equation}
        \begin{split}
            \displaystyle
            & H_{\zeta}^{}(\omega, \chi) = H_{\zeta c}^{}(\omega, \chi) + i H_{\zeta s}^{}(\omega, \chi)
            \\\\
            & \left\{\begin{split}
            \displaystyle
            & H_{\zeta c}^{}(\omega, \chi) = \dfrac{\dfrac{\omega^{2}L}{g} \cos\chi \sin\left(\dfrac{\omega^{2}L}{2g}\cos\chi\right)}{\pi^{2} - \left(\dfrac{\omega^{2}L}{2g}\cos\chi\right)^{2}}
            \\\\
            \displaystyle
            & H_{\zeta s}^{}(\omega, \chi) = 0
            \\\\
            \end{split}\right.
        \end{split}
    \end{equation}
    With the use of this transfer function, the spectrum of the effective wave can be obtained as follows:
    \begin{equation}
        \label{spectrum_effective_wave_theory}
        S_{\mathrm{eff}}^{}(\omega) = \left| H_{\zeta}^{}(\omega) \right|^{2} S_{\mathrm{w}}^{}(\omega)
    \end{equation}
    
    A more accurate approximation can be obtained by using a higher-order linear filter. In this study, the following 6th-order ARMA filter is used:
    \begin{equation}
        \label{ARMA6_all}
        \begin{split}
            \displaystyle &
            x_1^{(6)}+\alpha_1^{}x_1^{(5)}+\alpha_2^{}x_1^{(4)}+\alpha_3^{}x_1^{(3)}+\alpha_4^{}x_1^{(2)}+\alpha_5^{}x_1^{(1)}+\alpha_6^{}x_1^{(0)} 
            \\
            \displaystyle &
            =\sqrt{\pi}k\{w^{(1)}\}^{(3)} \,\,\,.
        \end{split}
    \end{equation}
    This expression is rewritten as:
     \begin{equation}
        \label{ARMA6}
        \left\{\begin{array}{l}
            \displaystyle
            \dfrac{\mathrm{d} x_{1}^{}}{\mathrm{d} t} = x_{2}^{} - \alpha_{1}^{} x_{1}^{}
            \\\\
            \displaystyle
            \dfrac{\mathrm{d} x_{2}^{}}{\mathrm{d} t} = x_{3}^{} - \alpha_{2}^{} x_{1}^{}
            \\\\
            \displaystyle
            \dfrac{\mathrm{d} x_{3}^{}}{\mathrm{d} t} = x_{4}^{} - \alpha_{3}^{} x_{1}^{} + \sqrt{\pi} k \dfrac{\mathrm{d} w(t)}{\mathrm{d} t}
            \\\\
            \displaystyle
            \dfrac{\mathrm{d} x_{4}^{}}{\mathrm{d} t} = x_{5}^{} - \alpha_{4}^{} x_{1}^{}
            \\\\
            \displaystyle
            \dfrac{\mathrm{d} x_{5}^{}}{\mathrm{d} t} = x_{6}^{} - \alpha_{5}^{} x_{1}^{}
            \\\\
            \displaystyle
            \dfrac{\mathrm{d} x_{6}^{}}{\mathrm{d} t} = - \alpha_{6}^{} x_{1}^{}
        \end{array}\right.
    \end{equation}
    The transfer function $H_{6}^{}(s)$ of the 6th-order linear filter is obtained from Eq.~(\ref{ARMA6}) as follows.
    \begin{equation}
        \label{ARMA6_transfer}
        H_{6}^{}(s)=\dfrac{\sqrt{\pi}ks^{3}}{s^{6}+s^{5}\alpha_{1}^{}+s^{4}\alpha_{2}^{}+s^{3}\alpha_{3}^{}+s^{2}\alpha_{4}^{}+s\alpha_{5}^{}+\alpha_{6}^{}}
    \end{equation}
    Finally, the 6th-order spectrum of the ARMA process can be obtained as:
    \begin{equation}
        \label{ARMA6_spectrum}
        \begin{split}
            \displaystyle &
            S_{6}(\omega)
            \\
            \displaystyle &
            =2\left| H_{6}(\omega) \right|^{2}\dfrac{1}{2\pi} =
            \\
            \displaystyle &
            \dfrac{k^{2}\omega^{6}}{\left(-\omega^{6} + \alpha_{2}^{}\omega^{4} - \alpha_{4}^{}\omega^{2} + \alpha_{6}^{} \right)^{2} + \left( \alpha_{1}^{}\omega^{5} - \alpha_{3}^{}\omega^{3} + \alpha_{5}^{}\omega^{}\right)^{2}}
        \end{split}
    .
    \end{equation}
    
    Next, it is necessary to obtain the coefficients of the ARMA process spectrum included in Eq.~(\ref{ARMA6_spectrum}) such that they well fit  the effective wave spectrum at a given $T_{01}$ and $H_{1/3}$ . It should be noted that even if these spectra have a good agreement, the system represented by Eq.~(\ref{ARMA6_spectrum}) can become unstable and problems in the modeling of time series may occur. Therefore, in this study, the stability criterion of the corresponding system is added as one of the conditions to determine the coefficients of the linear filter. The algebraic equation used for the stability criterion is given by Eq.~(\ref{ARMA6_algebra}).
    \begin{equation}
        \label{ARMA6_algebra}
        x^{6}+\alpha_1^{}x^{5}+\alpha_2^{}x^{4}+\alpha_3^{}x^{3}+\alpha_4^{}x^{2}+\alpha_5^{}x+\alpha_6^{} = 0
    \end{equation}
    If all the real parts of the poles of Eq.~(\ref{ARMA6_algebra}) are negative, then the system is stable.
    
\begin{figure}[h]
        \centering
        \includegraphics[scale=1]{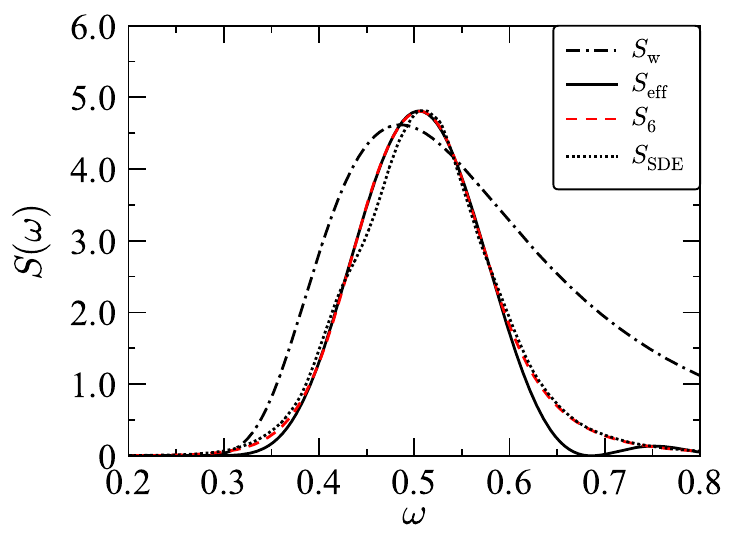}
        \caption{Comparison among ITTC spectrum : $S_{\mathrm{w}}^{}$, spectrum of effective wave represented by Eq.(\ref{spectrum_effective_wave_theory}) : $S_{\mathrm{eff}}^{}$, 6th-order ARMA spectrum : $S_{6}^{}$, and spectrum analysis result of time series obtained by solving SDE~(Eq.(\ref{ARMA6}) : $S_{\mathrm{SDE}}^{}$, sea state with $T_{01}=9.99[s]$ and $H_{1/3}=5.0[m]$ }
        \label{fig:spectrum_of_ARMA6_stability}
    \end{figure}

\begin{figure}[h]
        \centering
        \includegraphics[scale=0.95]{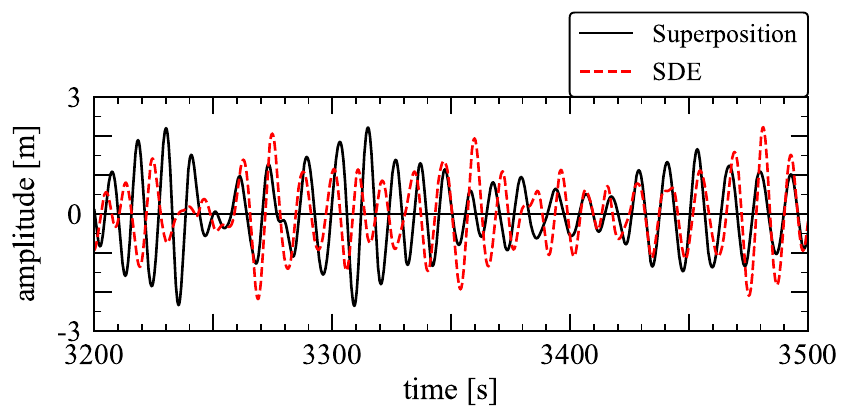}
        \caption{Comparison of time series of effective wave between numerical simulation by superposition methodology : Superposition and result obtained by solving SDE~(Eq.(\ref{ARMA6})) : SDE, sea state with $T_{01}=9.99[s]$ and $H_{1/3}=5.0[m]$}
        \label{fig:timeseries_of_ARMA6_stability}
    \end{figure}

    The computation results are shown in Table\ref{tab:ARMA_coefficient_stability} and Fig.~\ref{fig:spectrum_of_ARMA6_stability}. The solutions of Eq.~(\ref{ARMA6_algebra}) are obtained as {$-0.0861\pm0.432i$}, {$-0.237\pm0.422i$}, and {$-0.0909\pm0.547i$}. Since all real parts of these solutions are negative, it is clear that this system using the values from Table~\ref{tab:ARMA_coefficient_stability} is stable. In Fig.~\ref{fig:spectrum_of_ARMA6_stability}, it is observed that the 6th-order ARMA process spectrum shows a satisfactory agreement with the effective wave spectrum. Also, the time series of the effective wave is obtained by numerically solving the stochastic differential equation derived from Eq.~(\ref{ARMA6}) using the Euler-Maruyama's method. It can be seen from Fig.~\ref{fig:timeseries_of_ARMA6_stability} that the periods and amplitudes are almost the same between the result obtained from the superposition principle and the result obtained by solving the SDE~(Eq.(\ref{ARMA6})). In Fig.~\ref{fig:spectrum_of_ARMA6_stability}, the spectrum $S_{\mathrm{SDE}}^{}$ of the time series using the linear filter has a good agreement with $S_{6}^{}$. Therefore, it is clear that the time series of the effective wave can be modeled sufficiently by using the linear filter. The reason for a disagreement between random time series is that different random numbers are generated in each numerical calculation. However, the statistical properties, such as mean, variance mean frequency, and mean amplitude should coincide.

    \renewcommand{\arraystretch}{1.5}
    \begin{table}[h]
        \centering
        \caption{Coefficients of Eq.(\ref{ARMA6_spectrum})}
        \begin{tabular}{ll}
        \hline
         {\it Coefficient} & {\it ARMA6}
         \\
         \hline
          $\alpha_1^{}$ & $0.828$
          \\
          $\alpha_2^{}$ & $0.935$
          \\
          $\alpha_3^{}$ & $0.424$
          \\
          $\alpha_4^{}$ & $0.227$
          \\
          $\alpha_5^{}$ & $0.0490$
          \\
          $\alpha_6^{}$ & $0.0140$
          \\
          $k$ & $0.0459$
          \\
        \hline
        \end{tabular}
        \label{tab:ARMA_coefficient_stability}
    \end{table}
    \renewcommand{\arraystretch}{1.0}

    \begin{figure*}[h]
        \centering
        \includegraphics[scale=0.8]{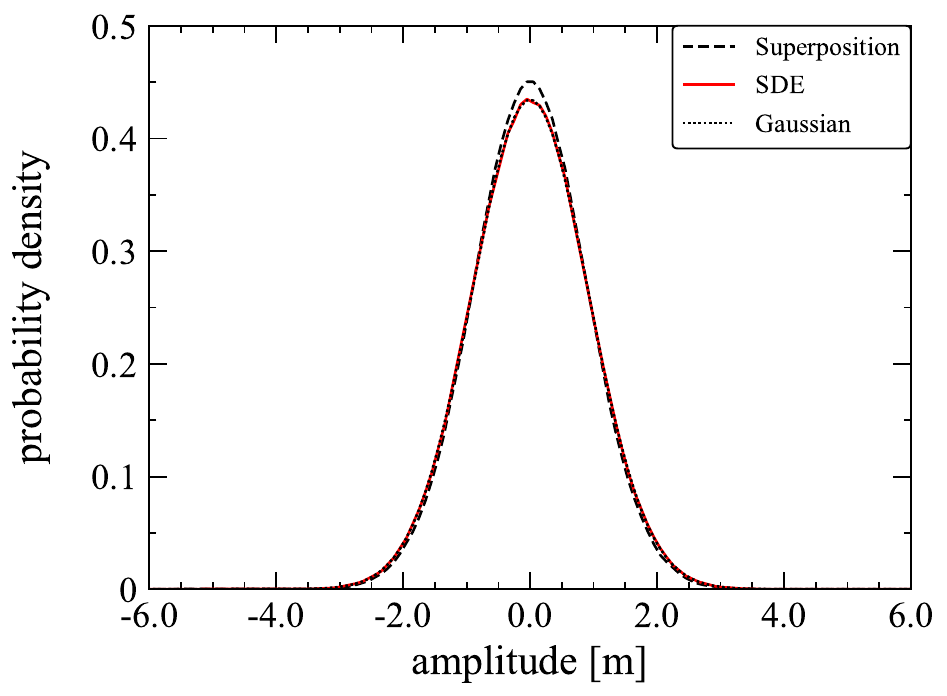}
        \includegraphics[scale=0.8]{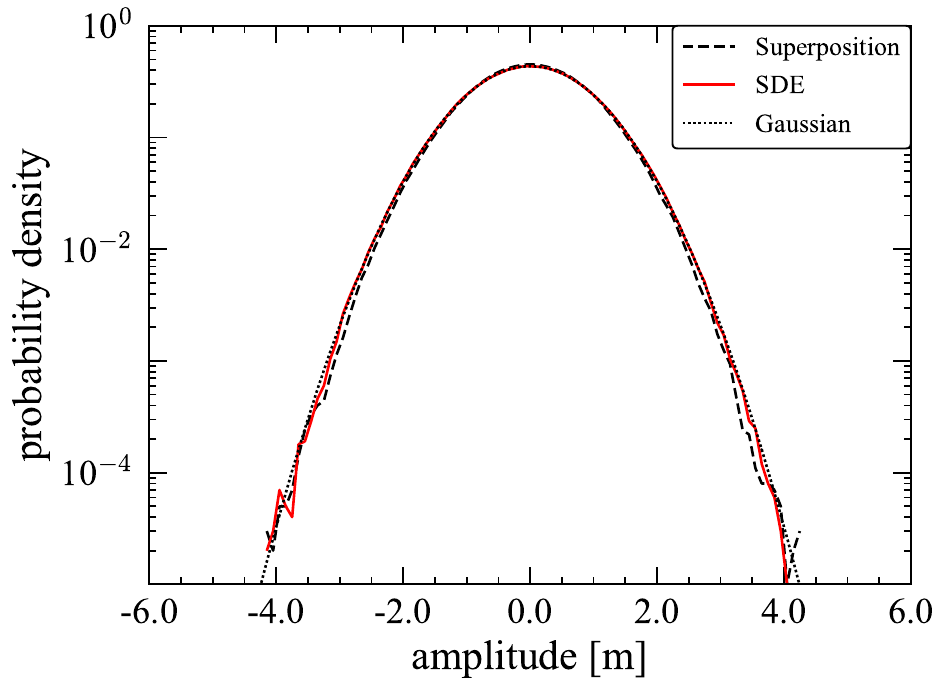}
        \caption{Comparison of PDF of effective wave amplitude among numerical simulation by superposition methodology : Superposition, result obtained by solving SDE Eq.(\ref{ARMA6}) : SDE, and Gaussian distribution by using moment equation result(Table\ref{tab:calculation_results_moment_value_and_superposition}) : Gaussian, sea state with $T_{01}=9.99[s]$ and $H_{1/3}=5.0[m]$}
        \label{fig:PDFofeffectivewave}
    \end{figure*}
    
    Next, in order to consider the GM variation in waves, the relationship between the amount of GM variation $\Delta$GM and wave amplitude at amidship is needed. The restoring arm for the case when the ship is heeling by two degrees in a regular wave, is calculated from hydrodynamic theory~\cite{Umeda1992} using a wavelength which is the same as the ship length. Then the wave crest or trough is set to be located at amidship, and GM is calculated for each wave amplitude. In Fig.~\ref{fig:HBL_vs_GM} it  the relationship between the amount of GM variation and wave amplitude at amidship can be seen. Here the wave amplitude is positive when the wave trough is located at amidship and the wave amplitude is negative when the wave crest is located at amidship. In this study, the actual data is approximated by means of a 12th-order polynomial. It can be observed in Fig.~\ref{fig:HBL_vs_GM} to have a satisfactory agreement.

    \begin{figure}[tb]
        \centering
        \includegraphics[scale=1]{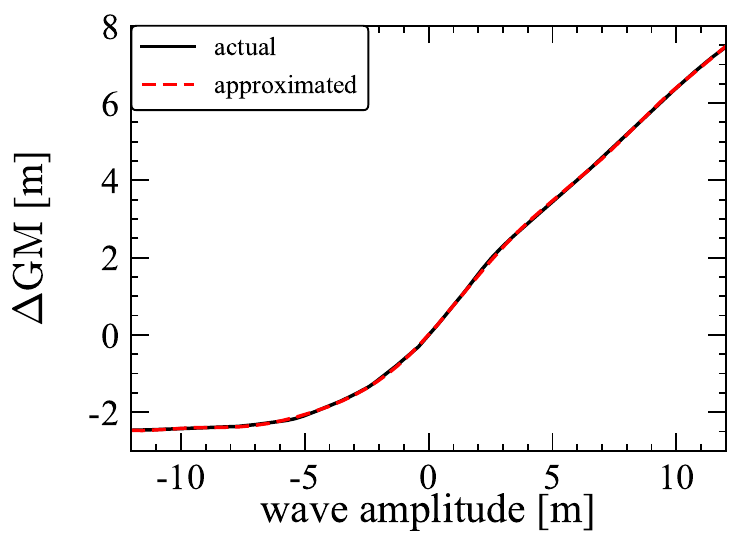}
        \caption{Relationship between $\Delta$GM and wave amplitude at amidship, subject ship : C11}
        \label{fig:HBL_vs_GM}
    \end{figure}
    \begin{figure}[tb]
        \centering
        \includegraphics[scale=0.95]{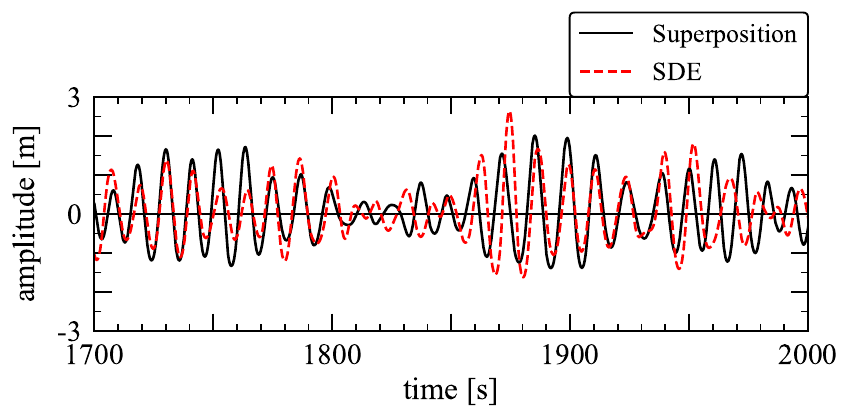}
        \caption{Comparison of time series of GM variation between numerical simulation by means of the superposition principle : Superposition and result obtained by solving SDE~(Eq.(\ref{ARMA6})) : SDE, sea state with $T_{01}=9.99[s]$ and $H_{1/3}=5.0[m]$}
        \label{fig:timeseriesofGM}
    \end{figure}
    
\begin{figure*}[tb]
        \centering
        \includegraphics[scale=0.8]{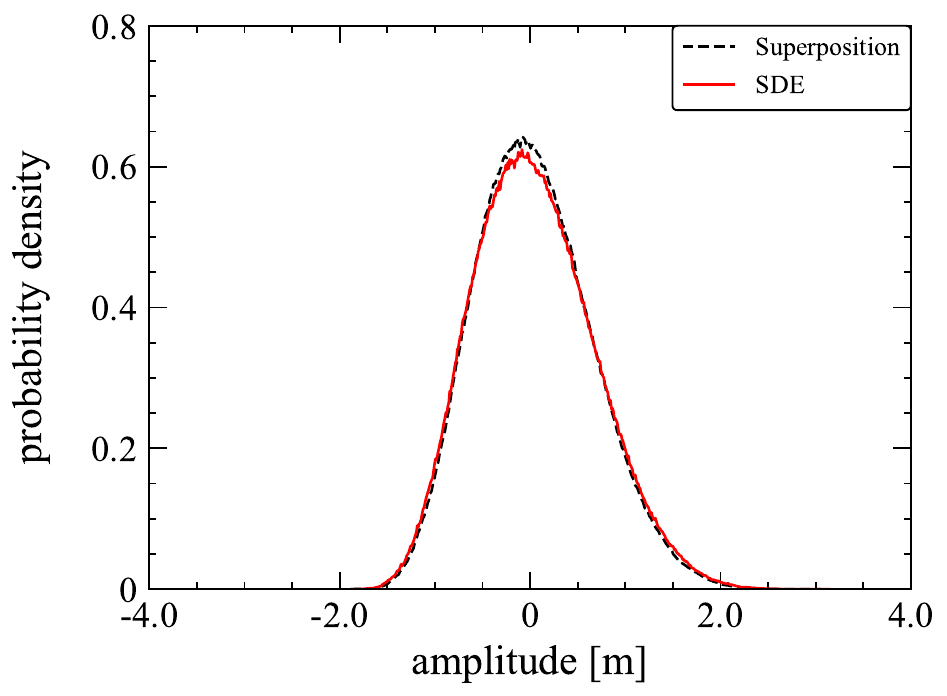}
        \includegraphics[scale=0.8]{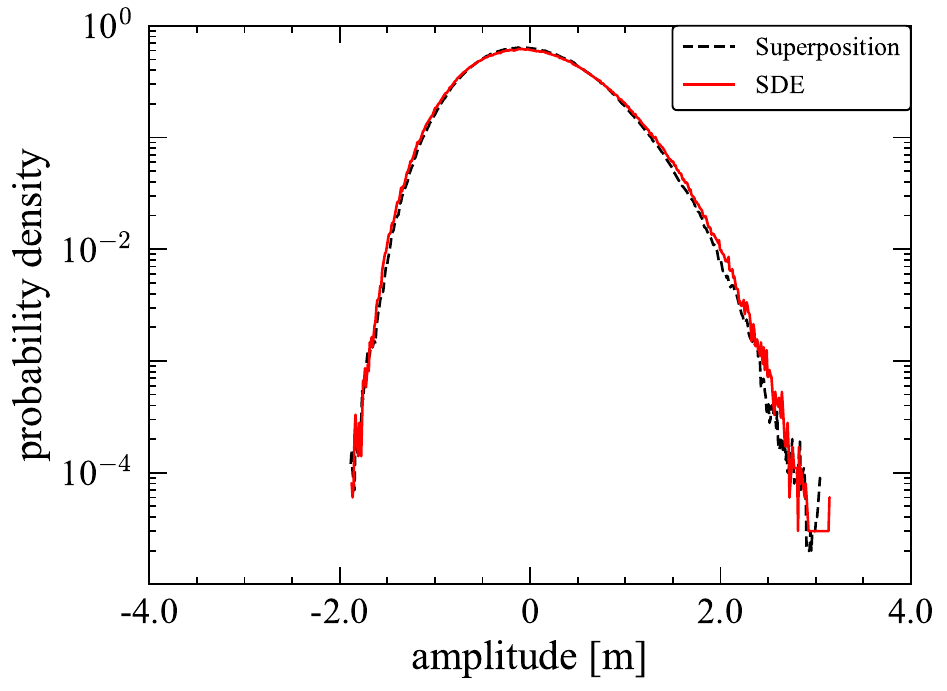}
        \caption{Comparison of PDF of GM variation between numerical simulation by means of the superposition principle : Superposition and result obtained by solving SDE Eq.~(\ref{ARMA6}) : SDE, sea state with $T_{01}=9.99[s]$ and $H_{1/3}=5.0[m]$}
        \label{fig:PDFofGM}
    \end{figure*}

\begin{figure}[h]
        \centering
        \includegraphics[scale=1]{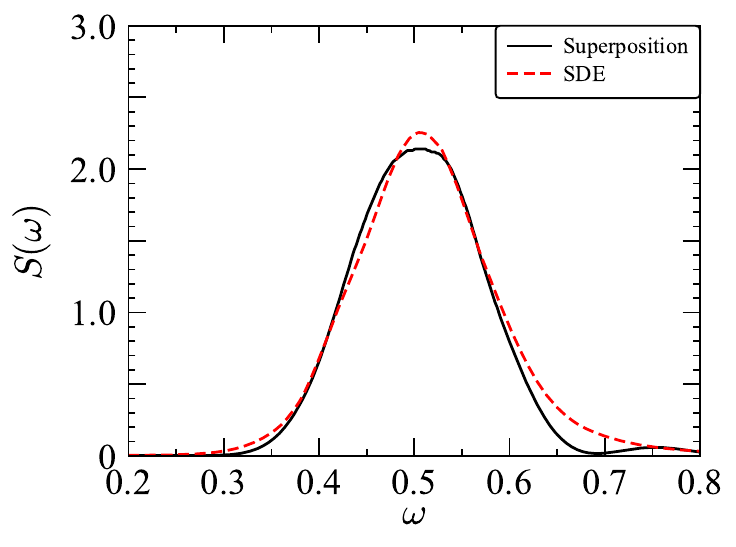}
        \caption{Comparison of spectrum analysis results using time series from numerical simulation by means of the superposition principle : Superposition and time series obtained by solving SDE~(Eq.(\ref{ARMA6}) and using the relationship in Fig.~\ref{fig:HBL_vs_GM}) : SDE, sea state with $T_{01}=9.99[s]$ and $H_{1/3}=5.0[m]$ }
        \label{fig:spectrum_of_GM}
    \end{figure}

    The time series of the GM variation is generated by using the time series of an effective wave and the relationship in Fig.~\ref{fig:HBL_vs_GM}. Fig.~\ref{fig:timeseriesofGM} shows an example of two time series of GM variation. In this figure, the black solid line shows the time series obtained by using the  wave superposition principle and the relationship in Fig.~\ref{fig:HBL_vs_GM}. The red solid line shows the time series from the result obtained by solving the stochastic differential equation derived from Eq.~(\ref{ARMA6}) and the relationship in Fig.~\ref{fig:HBL_vs_GM}. In this figure, it is observed that these periods and amplitudes are approximately equal. Therefore, it is clear that the time series of GM variation can be modeled sufficiently by using the linear filter and the relationship in Fig.~\ref{fig:HBL_vs_GM}. The reason for a disagreement of the time series of GM variation is that the time series of an effective wave obtained by each method is different. Also, Fig.~\ref{fig:PDFofGM} and \ref{fig:spectrum_of_GM} show satisfactory agreements about the probability density function of amplitude of GM variation and the spectrum of GM variation, which are obtained by each method, respectively.

    
\section{Moment Equation}
For phenomena induced by irregular external excitation, this process is modeled as Gaussian random process in many studies. This is because the distribution of this random process has generally a close shape to a Gaussian distribution and a statistical analysis becomes feasible. However, some of the phenomena induced by irregular excitation have a non-Gaussian distribution~\cite{YOKOYAMA2013}. The parametric rolling focused on in this study results from ship motion induced by irregular excitation. This phenomenon is a case which is non-Gaussian~\cite{Belenky2011}. 
The purpose of this paper is to determine the probability density function of roll angle, based on the moment values obtained by determining the moment equations. In order to obtain the moment equations, the system of the ship motion needs to be represented by a stochastic differential equation~(SDE)~\cite{Sobczyk1991}.

The procedure for deriving the moment equations from the stochastic differential equation is explained in the following. The It\^{o} stochastic differential equation is represented by 
    \begin{equation}
        \label{itoSDE}
        \mathrm{d} \bm{X}(t) = \bm{a}(t)\mathrm{d}t + \bm{b}(t)\mathrm{d}\bm{W}(t),
    \end{equation}
    where
    \begin{equation}
        \label{coe_itoSDE}
        \begin{array}{c}
            \displaystyle
            \bm{X}(t)=\left(
            \begin{array}{c}
                X_1(t) \\
                X_2(t) \\
                \vdots \\
                X_i(t) \\
            \end{array}
            \right)
            ,\,
            \bm{a}(t)=\left(
                \begin{array}{c}
                    a_1(t) \\
                    a_2(t) \\
                    \vdots \\
                    a_i(t) \\
                \end{array}
            \right)
            ,\\\\
            \displaystyle
            \bm{b}(t)=\left(
            \begin{array}{ccc}
                b_{11} & \cdots &  b_{1j} \\
                \vdots & \ddots & \vdots \\
                b_{i1} & \cdots &  b_{ij} \\
            \end{array}
            \right)
            ,\,
            \bm{W}(t)=\left(
            \begin{array}{c}
                W_1(t) \\
                W_2(t) \\
                \vdots \\
                W_j(t) \\
            \end{array}
            \right).
        \end{array}
    \end{equation}
    Here, $\bm{a}(t)$ is the drift term and $\bm{b}(t)$ is the diffusion term. Moreover, $\mathrm{d}\bm{W}(t)$ are independent increments of the standard Wiener process, whereby $\mathrm{d}\bm{W}(t)=\bm{W}(t+dt)-\bm{W}(t)$. When $f(t,\bm{X}(t))$ denotes a continuous function which is one-time differentiable with respect to $t$ and twice differentiable with respect to $\bm{X}(t)$, then the random process $f(t,\bm{X}(t))$ can be represented by using the It\^{o} formula as shown in the following equation.
    \begin{equation}
        \label{itofomulaSDE}
        \begin{array}{l}
            \displaystyle
            \mathrm{d} f(t,\bm{X}(t)) = 
            \left[
            \frac{\partial f}{\partial t}(t,\bm{X}(t)
            + \sum_{k=1}^{i}\frac{\partial}{\partial X_k}f(t,\bm{X}(t))a_k(t) \right.
            \\\\
            \displaystyle
            \left. + \frac{1}{2}\sum_{m=1}^{j}\sum_{k,l=1}^{i}\frac{\partial^2}{\partial X_k \partial X_l} f(t,\bm{X}(t))b_{km}(t)b_{lm}(t)
            \right] dt 
            \\\\
            \displaystyle
            +\sum_{m=1}^{j}\sum_{k=1}^{i}\frac{\partial f}{\partial X_k} (t,\bm{X}(t))b_{km}(t)\mathrm{d}W_m(t)
        \end{array}
    \end{equation}
    Taking the average on both sides of Eq.~(\ref{itofomulaSDE}) results in
    \begin{equation}
        \label{averagingitofomulaSDE}
        \begin{array}{l}
            \displaystyle
            \frac{\mathrm{d}}{\mathrm{d}t}\mathbb{E} \left[ f(\bm{X}) \right] =
            \sum_{k=1}^{i} \mathbb{E} \left[\frac{\partial f(\bm{X})}{\partial X_k} a_k(t) \right]
            \\\\
            \displaystyle
            + \frac{1}{2}\sum_{m=1}^{j}\sum_{k,l=1}^{i}\mathbb{E} \left[ \frac{ \partial^2 f(\bm{X}) }{ \partial X_k \partial X_l } b_{km}(t)b_{lm}(t) \right].
        \end{array}
    \end{equation}
    Here, $\mathbb{E}[\mathrm{d}W]=0$ is used, and $f(\bm{X})$ is denoted as 
    \begin{equation}
        \displaystyle f(\bm{X}) = \prod_{k=1}^{i} X_k^{C_k} = X_1^{C_1}X_2^{C_2} \cdots X_i^{C_i} . 
    \end{equation}
    Then, ordinary differential equations for the moments can be derived from Eq.~(\ref{averagingitofomulaSDE}).
    
    Next, the moment equations are derived from the equations of motion of parametric rolling~\cite{Maruyama2021}. In this study, the resulting system of equations is represented by the following 8th-order It\^o stochastic differential equation, which consists of a second-order SDE for the ship motion and a 6th-order SDE for the effective wave. 
    \begin{equation}
        \label{itopara}
        \left\{
        \begin{array}{l}
        \displaystyle \mathrm{d} X_1^{} = X_2^{} \mathrm{d} t \\
        \\
        \displaystyle \mathrm{d} X_2^{} = \left\{ - G(X_1^{},X_2^{}) - F(y_1^{}) X_1^{} \right\} \mathrm{d} t\\
        \\
        \displaystyle \mathrm{d} y_1^{} = \left( y_2^{} - \alpha_1^{} y_1^{} \right) \mathrm{d} t\\
        \\
        \displaystyle \mathrm{d} y_2^{} = \left( y_3^{} - \alpha_2^{} y_1^{} \right) \mathrm{d} t \\
        \\
        \displaystyle \mathrm{d} y_3^{} = \left( y_4^{} - \alpha_3^{} y_1^{} \right) \mathrm{d} t + \sqrt{\pi} k \mathrm{d} W(t), \\
        \\
        \displaystyle \mathrm{d} y_4^{} = \left( y_5^{} - \alpha_4^{} y_1^{} \right) \mathrm{d} t \\
        \\
        \displaystyle \mathrm{d} y_5^{} = \left( y_6^{} - \alpha_5^{} y_1^{} \right) \mathrm{d} t \\
        \\
        \displaystyle \mathrm{d} y_6^{} = - \alpha_6^{} y_1^{} \mathrm{d} t
        \end{array}
        \right.
    \end{equation}
    where
    \begin{equation}
        \begin{array}{l}
            \displaystyle
            G(X_1^{},X_2^{}) = \beta_1^{} X_2^{} + \beta_3^{} X_2^3 + \omega_{0}^{2} \sum_{n=1}^{5} \gamma_{2n-1} X_1^{2n-1},
            \\\\
            \displaystyle
            F(X_3^{}) = \frac{\omega_{0}^{2}}{\mathrm{GM}} \sum_{n=1}^{12} \rho_{n}^{} X_3^{n}.
        \end{array}
    \end{equation}
    Here $X_{1}^{}$, $X_{2}^{}$, and $X_{3}^{}$ are roll angle, roll velocity, and amplitude of the effective wave, respectively. Moreover, $\beta_{1}^{}$ is the linear and $\beta_{3}^{}$ is the cubic damping coefficient, divided by $I_{xx}^{}$, $I_{xx}^{}$ is the moment of inertia in roll (including the corresponding added moment of inertia), $\gamma_{i}^{}\,(i=1,3,5,7,and,9)$ is the coefficient of the i-th component of the polynomial fitted GZ curve, and $\rho_{i}^{}\,(i=1,2,\cdots,12)$ is the coefficient of the i-th component of the polynomial fitted the relationship between $\Delta$GM and wave height at midship in Fig.1.
    The notation of variables of the linear filter is rewritten, as $y_1=X_3,\, y_2=X_4,\, y_3=X_5,\, y_4=X_6,\, y_5=X_7,\, y_6=X_8$. Summarizing Eq.~(\ref{itopara}) in a matrix form like Eq.~(\ref{itoSDE}), the components of drift and diffusion are as follows:
    \begin{equation}
        \label{vecteritopara}
        \begin{array}{l}
            \displaystyle
            a_1^{} = X_2^{},\,a_2^{} = -G(X_1^{},X_2^{}) - F(X_3^{})X_1^{},\,a_3^{} = X_4^{} - \alpha_1^{} X_3^{},
            \\\\
            \displaystyle
            a_4^{} = X_5^{} - \alpha_2^{} X_3^{},\,a_5^{} = X_6^{} - \alpha_3^{} X_3^{},\,a_6^{} = X_7^{} - \alpha_4^{} X_3^{},
            \\\\
            \displaystyle
            a_7^{} = X_8^{} - \alpha_5^{} X_3^{},\,a_8^{} = -\alpha_6^{} X_3^{},
            \\\\
            \displaystyle
            b_{ij}^{} = 
            \begin{cases}
                \sqrt{\pi}k & (i=j=5)
                \\
                0      & \text{otherwise}
            \end{cases}
        \end{array}
        .
    \end{equation}
    By using $ f(\bm{X}) = \prod_{k=1}^{8} X_k^{C_k}$, 
    \begin{equation}
        \frac{\partial f(\bm{X})}{\partial X_k} = \frac{C_k}{X_k}f(\bm{X}),\,\,
        \frac{\partial^2 f(\bm{X})}{\partial X_5^2} = \frac{C_5(C_5-1)}{X_5^2}f(\bm{X})
        .
    \end{equation}
    Finally, the moment equation for the system in Eq.~(\ref{itopara}) is represented as:
    \begin{equation}
        \label{moment_para}
        \begin{array}{l}
            \displaystyle
            \frac{\mathrm{d}}{\mathrm{d}t}\mathbb{E} \left[ f(\bm{X}) \right] =
            \sum_{k=1}^{8} \mathbb{E} \left[\frac{\partial f(\bm{X})}{\partial X_k^{}} a_k^{}(t) \right]
            \\\\
            \displaystyle
            + \frac{1}{2}\sum_{m=1}^{8}\sum_{k,l=1}^{8}\mathbb{E} \left[ \frac{ \partial^2 f(\bm{X}) }{ \partial X_k^{} \partial X_l^{} } b_{km}^{}(t)b_{lm}^{}(t) \right]
            \\\\
            \displaystyle
            = C_1^{} \mathbb{E}\left[ X_1^{C_1^{}-1}X_2^{C_2^{}+1}X_3^{C_3^{}}X_4^{C_4^{}}X_5^{C_5^{}}X_6^{C_6^{}}X_7^{C_7^{}}X_8^{C_8^{}} \right]
            \\\\
            \displaystyle
            - C_2^{}\mathbb{E}\left[ G(X_1^{},X_2^{})X_1^{C_1^{}}X_2^{C_2^{}-1}X_3^{C_3^{}}X_4^{C_4^{}}X_5^{C_5^{}}X_6^{C_6^{}}X_7^{C_7^{}}X_8^{C_8^{}} \right]
            \\\\
            \displaystyle
            - C_2^{}\mathbb{E}\left[ F(X_3^{})X_1^{C_1^{}+1}X_2^{C_2^{}-1}X_3^{C_3^{}}X_4^{C_4^{}}X_5^{C_5^{}}X_6^{C_6^{}}X_7^{C_7^{}}X_8^{C_8^{}} \right]
            \\\\
            \displaystyle
            + C_3^{} \mathbb{E}\left[\left(X_4^{}-\alpha_1^{} X_3^{}\right)X_1^{C_1^{}}X_2^{C_2^{}}X_3^{C_3^{}-1}X_4^{C_4^{}}X_5^{C_5^{}}X_6^{C_6^{}}X_7^{C_7^{}}X_8^{C_8^{}}\right]
            \\\\
            \displaystyle
            + C_4^{} \mathbb{E}\left[\left(X_5^{}-\alpha_2^{} X_3^{}\right)X_1^{C_1^{}}X_2^{C_2^{}}X_3^{C_3^{}}X_4^{C_4^{}-1}X_5^{C_5^{}}X_6^{C_6^{}}X_7^{C_7^{}}X_8^{C_8^{}}\right]
            \\\\
            \displaystyle
            + C_5^{} \mathbb{E}\left[\left(X_6^{}-\alpha_3^{} X_3^{}\right)X_1^{C_1^{}}X_2^{C_2^{}}X_3^{C_3^{}}X_4^{C_4^{}}X_5^{C_5^{}-1}X_6^{C_6^{}}X_7^{C_7^{}}X_8^{C_8^{}}\right]
            \\\\
            \displaystyle
            + C_6^{} \mathbb{E}\left[\left(X_7^{}-\alpha_4^{} X_3^{}\right)X_1^{C_1^{}}X_2^{C_2^{}}X_3^{C_3^{}}X_4^{C_4^{}}X_5^{C_5^{}}X_6^{C_6^{}-1}X_7^{C_7^{}}X_8^{C_8^{}}\right]
            \\\\
            \displaystyle
            + C_7^{} \mathbb{E}\left[\left(X_8^{}-\alpha_5^{} X_3^{}\right)X_1^{C_1^{}}X_2^{C_2^{}}X_3^{C_3^{}}X_4^{C_4^{}}X_5^{C_5^{}}X_6^{C_6^{}}X_7^{C_7^{}-1}X_8^{C_8^{}}\right]
            \\\\
            \displaystyle
            - \alpha_6^{} C_8^{} \mathbb{E}\left[X_1^{C_1^{}}X_2^{C_2^{}}X_3^{C_3^{}+1}X_4^{C_4^{}}X_5^{C_5^{}}X_6^{C_6^{}}X_7^{C_7^{}}X_8^{C_8^{}-1}\right]
            \\\\
            \displaystyle
            + \frac{\pi k^2}{2} C_5^{}(C_5^{}-1) \mathbb{E}\left[ X_1^{C_1^{}}X_2^{C_2^{}}X_3^{C_3^{}}X_4^{C_4^{}}X_5^{C_5^{}-2}X_6^{C_6^{}}X_7^{C_7^{}}X_8^{C_8^{}} \right]
        \end{array}
    \end{equation}
    From Eq.~(\ref{moment_para}), the following eight first order moment equations are obtained:
    \begin{equation}
        \label{first_moment_para}
        \begin{array}{l}
            \displaystyle
            \frac{\mathrm{d}}{\mathrm{d}t}\mathbb{E}\left[X_1^{}\right] = \mathbb{E}\left[X_2^{}\right]
            \\\\
            \displaystyle
            \frac{\mathrm{d}}{\mathrm{d}t}\mathbb{E}\left[X_2^{}\right] = - \mathbb{E}\left[G(X_1^{},X_2^{})\right] - \mathbb{E}\left[F(X_3^{})X_1^{}\right]
            \\\\
            \displaystyle
            \frac{\mathrm{d}}{\mathrm{d}t}\mathbb{E}\left[X_i^{}\right] = \mathbb{E}\left[X_{i+1}^{}\right] - \alpha_{i-2}^{} \mathbb{E}\left[X_3^{}\right]\,\,\,(i=3\cdots7)
            \\\\
            \displaystyle
            \frac{\mathrm{d}}{\mathrm{d}t}\mathbb{E}\left[X_8^{}\right] = - \alpha_6^{} \mathbb{E}\left[X_3^{}\right]
        \end{array}
    \end{equation}
    Moreover, thirty-six second order moment equations are obtained as: 
    \begin{equation*}
        \begin{array}{l}
            \displaystyle
            \frac{\mathrm{d}}{\mathrm{d}t}\mathbb{E}\left[X_1^2\right] = 2 \mathbb{E}\left[X_1^{}X_2^{}\right]
            \\\\
            \displaystyle
            \frac{\mathrm{d}}{\mathrm{d}t}\mathbb{E}\left[X_1^{}X_2^{}\right] = \mathbb{E}\left[X_2^2\right] - \mathbb{E}\left[G(X_1^{},X_2^{})X_1^{}\right] - \mathbb{E}\left[F(X_3^{})X_1^2\right]
            \\\\
            \displaystyle
            \frac{\mathrm{d}}{\mathrm{d}t}\mathbb{E}\left[X_1^{}X_i^{}\right] = \mathbb{E}\left[X_2^{}X_i^{}\right] + \mathbb{E}\left[X_1^{}X_{i+1}^{}\right]
            \\
            \displaystyle
            \quad\quad\quad\quad\quad\quad
            - \alpha_{i-2}^{} \mathbb{E}\left[X_1^{}X_3^{}\right]\,\,\,(i=3\cdots7)
            \\\\
            \displaystyle
            \frac{\mathrm{d}}{\mathrm{d}t}\mathbb{E}\left[X_1^{}X_8^{}\right] = \mathbb{E}\left[X_2^{}X_8^{}\right] - \alpha_6^{} \mathbb{E}\left[X_1^{}X_3^{}\right]
            \\\\
            \displaystyle
            \frac{\mathrm{d}}{\mathrm{d}t}\mathbb{E}\left[X_2^2\right] = - 2 \mathbb{E}\left[G(X_1^{},X_2^{})X_2^{}\right] - 2 \mathbb{E}\left[F(X_3^{})X_1^{}X_2^{}\right]
            \\\\
            \displaystyle
            \frac{\mathrm{d}}{\mathrm{d}t}\mathbb{E}\left[X_2^{}X_i^{}\right] = - \mathbb{E}\left[G(X_1^{},X_2^{})X_i^{}\right] - \mathbb{E}\left[F(X_3^{})X_1^{}X_i^{}\right] 
            \\
            \displaystyle
            \quad\quad\quad\quad\quad\quad
            + \mathbb{E}\left[X_2^{}X_{i+1}^{}\right] - \alpha_{i-2}^{} \mathbb{E}\left[X_2^{}X_3^{}\right]\,\,\,(i=3\cdots7)
            \\\\
            \displaystyle
            \frac{\mathrm{d}}{\mathrm{d}t}\mathbb{E}\left[X_2^{}X_8^{}\right] = - \mathbb{E}\left[G(X_1^{},X_2^{})X_8^{}\right] - \mathbb{E}\left[F(X_3^{})X_1^{}X_8^{}\right]
            \\
            \displaystyle
            \quad\quad\quad\quad\quad\quad
            - \alpha_6^{} \mathbb{E}\left[X_2^{}X_3^{}\right]
            \\\\
            \displaystyle
            \frac{\mathrm{d}}{\mathrm{d}t}\mathbb{E}\left[X_3^2\right] = 2\mathbb{E}\left[X_3^{}X_4^{}\right] - 2\alpha_1^{}\mathbb{E}\left[X_3^2\right]
            \\\\
            \displaystyle
            \frac{\mathrm{d}}{\mathrm{d}t}\mathbb{E}\left[X_3^{}X_{i}^{}\right] = \mathbb{E}\left[X_4^{}X_{i}^{}\right] - \alpha_1^{}\mathbb{E}\left[X_3^{}X_{i}^{}\right] 
            \\
            \displaystyle
            \quad\quad\quad\quad\quad\quad
            + \mathbb{E}\left[X_3^{}X_{i+1}^{}\right] - \alpha_{i-2}^{}\mathbb{E}\left[X_3^2\right]\,\,\,(i=4\cdots7)
            \\\\
            \displaystyle
            \frac{\mathrm{d}}{\mathrm{d}t}\mathbb{E}\left[X_3^{}X_8^{}\right] = \mathbb{E}\left[X_4^{}X_8^{}\right] - \alpha_1^{}\mathbb{E}\left[X_3^{}X_8^{}\right] - \alpha_6^{}\mathbb{E}\left[X_3^2\right]
            \\\\
            \displaystyle
            \frac{\mathrm{d}}{\mathrm{d}t}\mathbb{E}\left[X_4^2\right] = 2\mathbb{E}\left[X_4^{}X_5^{}\right] - 2\alpha_2^{}\mathbb{E}\left[X_3^{}X_4^{}\right]
            \\\\
            \displaystyle
            \frac{\mathrm{d}}{\mathrm{d}t}\mathbb{E}\left[X_4^{}X_i^{}\right] = \mathbb{E}\left[X_5^{}X_i^{}\right] - \alpha_2^{}\mathbb{E}\left[X_3^{}X_i^{}\right]
            \\
            \displaystyle
            \quad\quad\quad\quad\quad\quad
            + \mathbb{E}\left[X_4^{}X_{i+1}^{}\right] - \alpha_{i-2}^{}\mathbb{E}\left[X_3^{}X_4^{}\right]\,\,\,(i=5\cdots7)
            \\\\
            \displaystyle
            \frac{\mathrm{d}}{\mathrm{d}t}\mathbb{E}\left[X_4^{}X_8^{}\right] = \mathbb{E}\left[X_5^{}X_8^{}\right] - \alpha_2^{}\mathbb{E}\left[X_3^{}X_8^{}\right] - \alpha_6^{}\mathbb{E}\left[X_3^{}X_4^{}\right]
        \end{array}
    \end{equation*}
    \begin{equation}
        \label{second_moment_para}
        \begin{array}{l}
            \displaystyle
            \frac{\mathrm{d}}{\mathrm{d}t}\mathbb{E}\left[X_5^2\right] = 2\mathbb{E}\left[X_5^{}X_6^{}\right] - 2\alpha_3^{}\mathbb{E}\left[X_3^{}X_5^{}\right] + \pi k_{}^{2}
            \\\\
            \displaystyle
            \frac{\mathrm{d}}{\mathrm{d}t}\mathbb{E}\left[X_5^{}X_i^{}\right] = \mathbb{E}\left[X_6^{}X_i^{}\right] - \alpha_3^{}\mathbb{E}\left[X_3^{}X_i^{}\right]
            \\
            \displaystyle
            \quad\quad\quad\quad\quad\quad
            + \mathbb{E}\left[X_5^{}X_{i+1}^{}\right] - \alpha_{i-2}^{}\mathbb{E}\left[X_3^{}X_5^{}\right]\,\,\,(i=6,7)
            \\\\
            \displaystyle
            \frac{\mathrm{d}}{\mathrm{d}t}\mathbb{E}\left[X_5^{}X_8^{}\right] = \mathbb{E}\left[X_6^{}X_8^{}\right] - \alpha_3^{}\mathbb{E}\left[X_3^{}X_8^{}\right] - \alpha_6^{}\mathbb{E}\left[X_3^{}X_5^{}\right]
            \\\\
            \displaystyle
            \frac{\mathrm{d}}{\mathrm{d}t}\mathbb{E}\left[X_6^2\right] = 2\mathbb{E}\left[X_6^{}X_7^{}\right] - 2\alpha_4^{}\mathbb{E}\left[X_3^{}X_6^{}\right]
            \\\\
            \displaystyle
            \frac{\mathrm{d}}{\mathrm{d}t}\mathbb{E}\left[X_6^{}X_7^{}\right] = \mathbb{E}\left[X_7^2\right] - \alpha_4^{}\mathbb{E}\left[X_3^{}X_7^{}\right]
            \\
            \displaystyle
            \quad\quad\quad\quad\quad\quad
            + \mathbb{E}\left[X_6^{}X_8^{}\right] - \alpha_5^{}\mathbb{E}\left[X_3^{}X_6^{}\right]
            \\\\
            \displaystyle
            \frac{\mathrm{d}}{\mathrm{d}t}\mathbb{E}\left[X_6^{}X_8^{}\right] = \mathbb{E}\left[X_7^{}X_8^{}\right] - \alpha_4^{}\mathbb{E}\left[X_3^{}X_8^{}\right] - \alpha_6^{}\mathbb{E}\left[X_3^{}X_6^{}\right]
            \\\\
            \displaystyle
            \frac{\mathrm{d}}{\mathrm{d}t}\mathbb{E}\left[X_7^2\right] = 2\mathbb{E}\left[X_7^{}X_8^{}\right] - 2\alpha_5^{}\mathbb{E}\left[X_3^{}X_7^{}\right]
            \\\\
            \displaystyle
            \frac{\mathrm{d}}{\mathrm{d}t}\mathbb{E}\left[X_7^{}X_8^{}\right] = \mathbb{E}\left[X_8^2\right] - \alpha_5^{}\mathbb{E}\left[X_3^{}X_8^{}\right] - \alpha_6^{}\mathbb{E}\left[X_3^{}X_7^{}\right]
            \\\\
            \displaystyle
            \frac{\mathrm{d}}{\mathrm{d}t}\mathbb{E}\left[X_8^2\right] = - 2\alpha_6^{}\mathbb{E}\left[X_3^{}X_8^{}\right]
        \end{array}
    \end{equation}
In general, a nonlinear system generates an infinite hierarchy of moment equations. In order to form a closed set of moment equations, higher-order moments need to be truncated. Therefore, in this study the cumulant neglect closure method\cite{Sun1987}\cite{Sun1989}\cite{Wojtkiewicz1996} is used. In particular, the second-order cumulant neglect method is used, which ignores cumulants of higher than the third-orde. This is closure method is the same as the Gaussian closure method, because the third-order and higher order cumulants of a Gaussian distribution are zero. In this study, the first and second-order moment equations include third and higher-order moments. Hence, the third and higher-order moments need to be represented using first and second-order moments. Comparing the coefficients of series expansions of a moment generating function and a cumulant generating function, the relations between moments and cumulants can be obtained, which is shown in detail in Appendix~1.

\section{Calculation results}
In this study, the numerical simulations by means of the superposition principle and by solving the SDE are analysed. And the moment equations are solved. Furthermore, two types of the arbitrary PDF and the objective function are set.

\subsection{Subject Ship}\label{sec:sec4sub1}
\begin{figure}[h]
    \centering
    \includegraphics[scale=1]{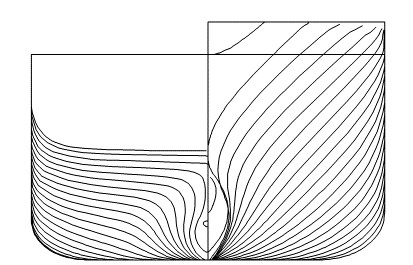}
    \caption{Body Plan (C11)}
    \label{fig:bodyplanC11}
\end{figure}

\begin{figure}[h]
    \centering
    \includegraphics[scale=0.8]{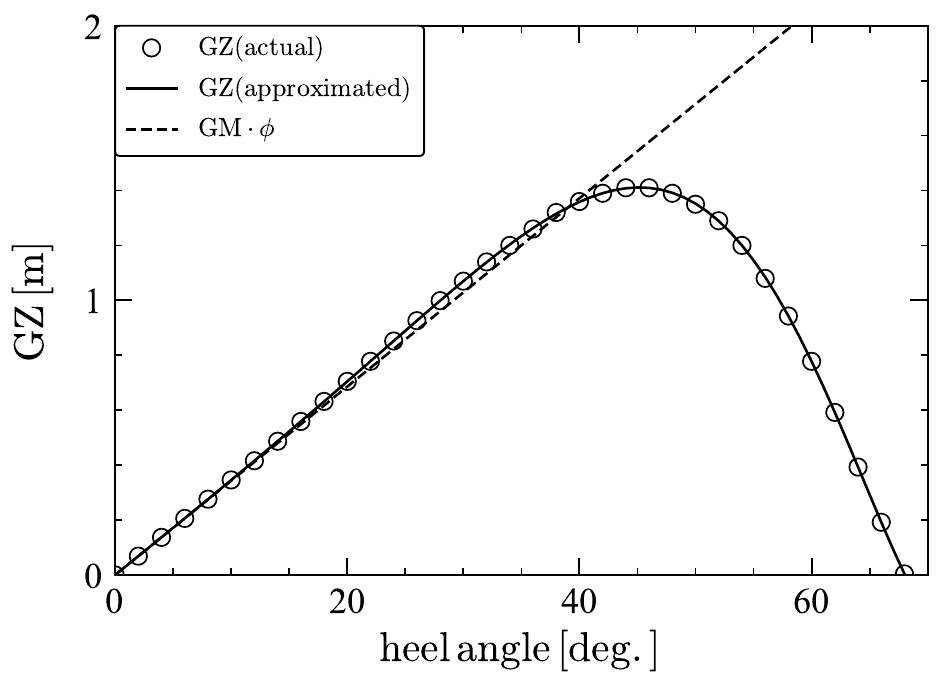}
    \caption{GZ curve in still water (C11)}
    \label{fig:GZcurve_stillwater_C11}
\end{figure}

In order to validate the moment equation method, the results are compared with numerical simulations. In this study, for a subject ship, the calculation results are compared for the Froude number Fn = 0.0 in head seas. Here, the C11 class post-Panamax containership~\cite{Hashimoto2010} is used as the subject ship. The body plan and principal particulars of this ship is shown in Fig.~\ref{fig:bodyplanC11} and Table~\ref{tab:principal_C11}, respectively. Fig.~\ref{fig:GZcurve_stillwater_C11} shows the GZ curve for the C11 containership, resulting from hydrostatic calculations in calm water. In this study, the actual GZ is approximated by a 9th order polynomial in order to provide a reasonable GZ curve. As shown in Fig.~\ref{fig:GZcurve_stillwater_C11}, the C11 containership has a GZ curve which is linear up to an roll angle of around 40 degrees. The calculation condition in this study is that the wave mean period $T_{01}^{}$ is $9.99$[s] and the significant wave height $H_{1/3}^{}$ is $
5.0$[m]. Also, the roll damping coefficients are estimated by Ikeda’s simplified method~\cite{Kawahara2012}. As a result, the damping coefficients are obtained as $\beta_1^{}=3.64\times10_{}^{-3}$ and $\beta_3^{}=4.25$.

\renewcommand{\arraystretch}{1.5}
\begin{table}[h]
  \begin{center}
    \caption{Principal Particulars of the subject ship at full scale}
    \begin{tabular}{ll}
    \hline
     {\it Items} & {\it C11}
     \\
     \hline
      Length:\,$L_{pp}$ & $262.0 \,[\mathrm{m}]$
      \\
      Breadth:\,$B$ & $40.0 \,[\mathrm{m}]$
      \\
      Depth:\,$D$ & $24.45 \,[\mathrm{m}]$
      \\
      Draught:\,$d$ & $11.5 \,[\mathrm{m}]$
      \\
      Block coefficient:\,$C_b$ & $0.562$
      \\
      Metacentric height:\,$GM$ & $1.965 \,[\mathrm{m}]$
      \\
      Natural roll period:\,$T_{\phi}$ & $25.1 \,[\mathrm{s}]$
      \\
      Bilge keel length ratio:\,$L_{BK}/L_{pp}$ & $0.292$
      \\
      Bilge keel breadth ratio:\,$B_{BK}/B$ & $0.0100$
      \\
    \hline
    \end{tabular}
    \label{tab:principal_C11}
  \end{center}
\end{table}
\renewcommand{\arraystretch}{1.0}

\subsection{Calculation of moment equations}\label{sec:sec4sub2}
In order to compare the moment values obtained by solving the moment equations, the numerical simulation of the superposition principle and the SDE are explained first. The numerical simulation using the superposition principle is computed by means of the 4th order Runge–Kutta method, while the numerical simulation of the SDE~Eq.~(\ref{itopara}) is computed using the Euler-Maruyama method. The time step is set such that the former is 0.02[s] and the latter is 0.001[s]. Both initial conditions are set to a roll angle of 5[deg.] and a roll velocity of 0[deg./s]. The number of realizations of the numerical simulation is 100, and each simulation time is 1 hour. 

\begin{figure}[h]
        \centering
        \includegraphics[scale=1]{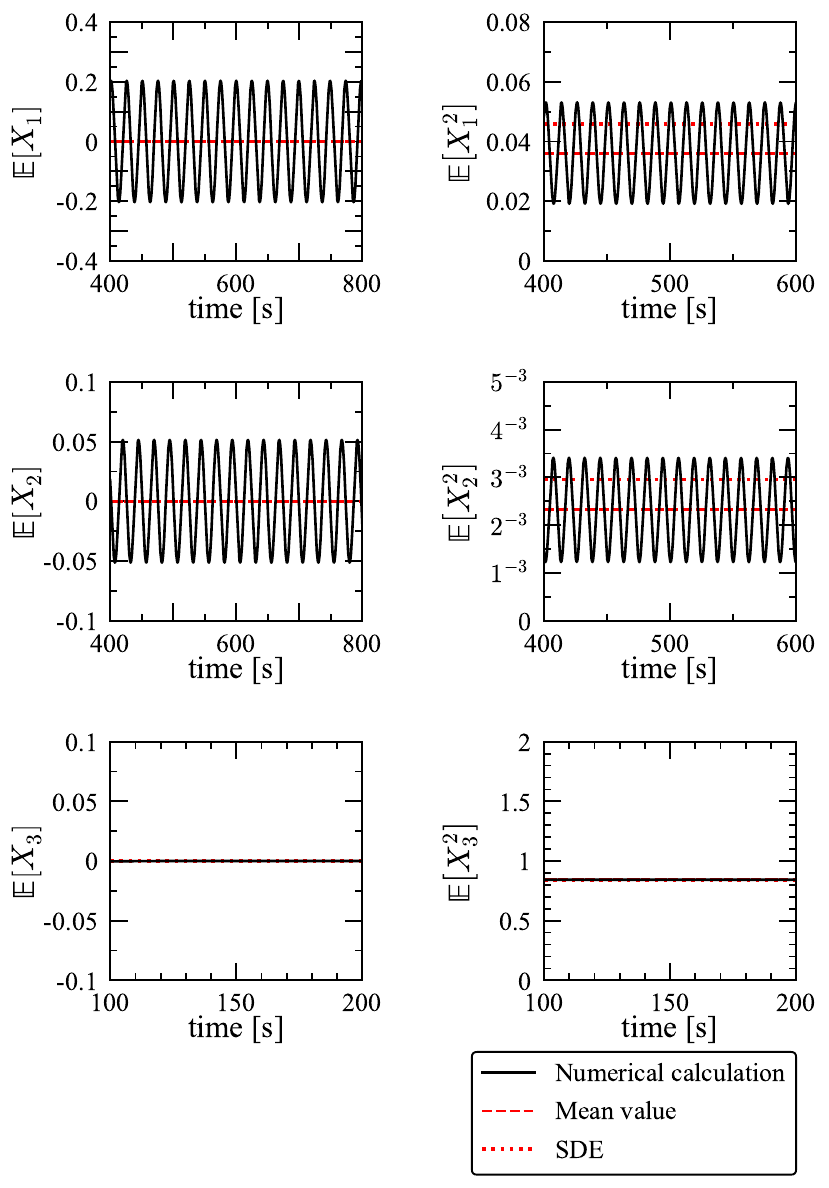}
        \caption{First and second order moment of $X_{1}^{},X_{2}^{},X_{3}^{}$ obtained by computing the ordinary differential equations for  moments (Eq.~(\ref{first_moment_para}) and (\ref{second_moment_para})), and moment values from SDE (Table~\ref{tab:calculation_results_superposition_and_SDE}) }
        \label{fig:moment_value_in_unstable}
    \end{figure}

Next, the ordinary differential equations for the moment equations ( Eq.~(\ref{first_moment_para}) and (\ref{second_moment_para}) ) are explained. These equations are computed using the 4th order Runge–Kutta method as an initial value problem. The initial condition for all moments is set to a value of 0.01. Here, the time step is set to 0.01[s]. Fig.~\ref{fig:moment_value_in_unstable} shows the computation results. The moment values in Table~\ref{tab:calculation_results_moment_value_and_superposition} are determined by taking the mean of the steady state time series of the respective moments. The red dashed line shows these mean values. The blue dashed line shows the moment values obtained by computing the SDE numerically.
 
In this study, it should be noted that the moment equations are calculated in an unsteady state. The calculation of moment equations in steady-state needs to solve simultaneous nonlinear equations. Then the solutions can be obtained by using the Newton-Raphson method and Jacobian matrix. The convergence and the corresponding matrix calculation are complex. On the other hand, the region of steady state can be determined easily from the computation of the corresponding ordinary differential equation. Furthermore, we consider that it is an appropriate method from the perspective of unaffected on the number of moment equations.
 
As shown in Fig.~\ref{fig:moment_value_in_unstable}, the solution in steady-state is oscillatory. When the cumulant neglect order is of third order, it is noteworthy that this oscillation vanishes and the moment values become close to the moment values derived from numerical computations of the SDE, as shown in Appendix~2.
Concerning the wave amplitude $X_3^{}$, it can be seen from Table~\ref{tab:calculation_results_moment_value_and_superposition} that the moment values obtained by solving the SDE part consisting of the linear filter are almost equal to the moment values obtained by solving the moment equations. On the other hand, comparing the moment values of $X_3^{}$ obtained by using the superposition principle with the moment values obtained by solving the linear filter SDE, the difference is about 7$\%$. We consider the cause of this to be caused by the discrepancy between the ARMA spectrum and the effective wave spectrum in the region from $\omega=0.6$ to $\omega=0.8$, which can be seen in Fig.~\ref{fig:spectrum_of_ARMA6_stability}. It is thus possible that these moment values are improved due to eliminating this discrepancy. Fig.~\ref{fig:PDFofeffectivewave} shows the comparison of the PDF of the effective wave amplitude. It is observed that the Gaussian distribution using mean and variance obtained from the moment equations has a good agreement with the PDF obtained by solving the linear filter SDE.

    \renewcommand{\arraystretch}{1.5}
     \begin{table}[h]
      \begin{center}
        \caption{Moment values obtained by solving SDE~(Eq.(\ref{itopara})) and numerical simulation using the superposition principle. }
        \begin{tabular}{lll}
        \hline
         {} &  {Superposition} &  {SDE}
         \\
         \hline
          $\mathbb{E}\left[X_{1}^{\,}\right]$ & $-4.01\times10_{}^{-5}$ & $4.99\times10_{}^{-6}$
          \\
          $\mathbb{E}\left[X_{1}^{2}\right]$ & $4.38\times10_{}^{-2}$ & $4.58\times10_{}^{-2}$
          \\
          $\mathbb{E}\left[X_{2}^{\,}\right]$ & $-6.00\times10_{}^{-6}$ & $-3.61\times10_{}^{-6}$
          \\
          $\mathbb{E}\left[X_{2}^{2}\right]$ & $2.81\times10_{}^{-3}$ & $2.95\times10_{}^{-3}$
          \\
          $\mathbb{E}\left[X_{3}^{\,}\right]$ & $-1.66\times10_{}^{-5}$ & $-7.63\times10_{}^{-5}$
          \\
          $\mathbb{E}\left[X_{3}^{2}\right]$ & $0.786$ & $0.842$
          \\
        \hline
        \end{tabular}
        \label{tab:calculation_results_superposition_and_SDE}
      \end{center}
    \end{table}
    \renewcommand{\arraystretch}{1.0}

    \renewcommand{\arraystretch}{1.5}
    \begin{table}[h]
      \begin{center}
        \caption{Moment values obtained by solving the moment equations (Eq.(\ref{first_moment_para}) and (\ref{second_moment_para}) ) and by solving SDE~(Eq.(\ref{itopara})) }
        \begin{tabular}{lll}
        \hline
         {} & {SDE} & {Moment equations}
         \\
         \hline
          $\mathbb{E}\left[X_{1}^{}\right]$ & $4.99\times10_{}^{-6}$  & $-1.02\times10_{}^{-5}$
          \\
          $\mathbb{E}\left[X_{1}^{2}\right]$ & $4.58\times10_{}^{-2}$ & $3.60\times10_{}^{-2}$
          \\
          $\mathbb{E}\left[X_{2}^{}\right]$ & $-3.61\times10_{}^{-6}$ & $-3.13\times10_{}^{-4}$
          \\
          $\mathbb{E}\left[X_{2}^{2}\right]$ & $2.95\times10_{}^{-3}$ & $2.32\times10_{}^{-3}$
          \\
          $\mathbb{E}\left[X_{3}^{}\right]$ & $-7.63\times10_{}^{-5}$ & $0.00$
          \\
          $\mathbb{E}\left[X_{3}^{2}\right]$ & $0.842$ & $0.843$
          \\
        \hline
        \end{tabular}
        \label{tab:calculation_results_moment_value_and_superposition}
      \end{center}
    \end{table}
    \renewcommand{\arraystretch}{1.0}
    
\subsection{Procedure for determining the PDF, and results}\label{sec:sec4sub3}
Based on the moment values in Table\ref{tab:calculation_results_moment_value_and_superposition} by computing the moment equations, the PDF of roll angle is determined. In this study, the following non-Gaussian PDF shape types are set:
    \begin{equation}
        \label{PDF_arbitrary_type1}
        \begin{split}
            \displaystyle
            &\text{type1 : }
            \\
            \displaystyle
            & \mathcal{P}_{1}^{}(X_1^{})
            \\
            \displaystyle
            & = C \exp{ \left\{ - \left( d_1^{} X_1^{} + d_2^{} X_{1}^{2} + d_3^{} X_{1}^{3} + d_4^{} X_{1}^{4} \right) \right\} }
        \end{split}
    \end{equation}
    \begin{equation}
        \label{PDF_arbitrary_type2}
        \begin{split}
            \displaystyle
            &\text{type2 : }
            \\
            \displaystyle
            & \mathcal{P}_{2}^{}(X_1^{}) 
            \\
            \displaystyle
            & = C \exp{ \left\{ - \left( d_1^{}\left| X_1^{} \right| + d_2^{} \left| X_1^{} \right|_{}^{2} + d_3^{} \left| X_1^{} \right|_{}^{3} + d_4^{} \left| X_1^{} \right|_{}^{4} \right) \right\} }
        \end{split}
    \end{equation}
    Here, $C$ is a normalization constant. In order to determine the coefficients of Eqs.~(\ref{PDF_arbitrary_type1}) and (\ref{PDF_arbitrary_type2}), the following expression is suggested. 
    \begin{equation}
        \label{Estimation_PDF_arbitrary1}
        J_{n}^{} = \int_{-\,\infty}^{+\infty} X_{1}^{n} \mathcal{P}(X)\,dX - \mathbb{E}\left[X_{1}^{n}\right]
    \end{equation}
    The moment values, which are obtained by solving the moment equations in Eq.~(\ref{first_moment_para}) and (\ref{second_moment_para}), are the first- and the second-order moments, only. 
    Hence, higher-order moment values need to be obtained by using the cumulant neglect closure method as shown in Appendix 1. Furthermore, the following objective function $J\left( d_{1}^{},d_{2}^{},d_{3}^{},d_{4}^{}\right)$ is set.
    \begin{equation}
        \label{Estimation_PDF_arbitrary2}
        J\left( d_{1}^{},d_{2}^{},d_{3}^{},d_{4}^{}\right) = \sum_{i=1}^{4} l_{i}^{}\,J_{i}^{}
    \end{equation}
    Here, $l_{i}^{}$ are weights and $l_{i}^{}=1\,(i = 1, 2,...,  4)$. 
    
\begin{figure}[h]
        \centering
        \includegraphics[scale=0.95]{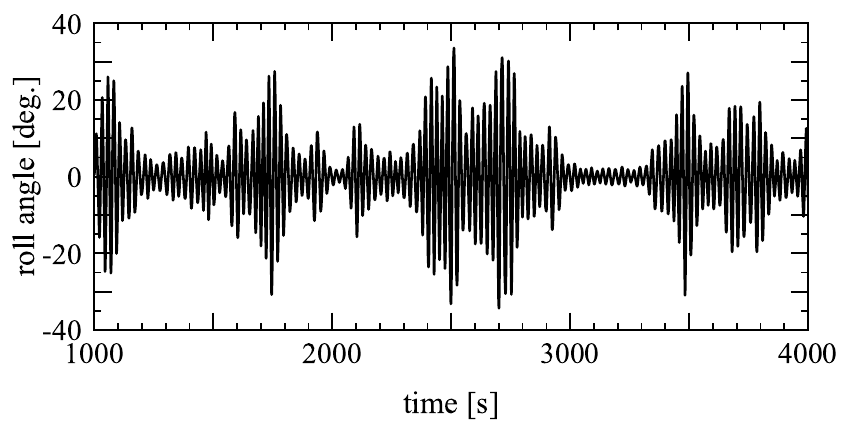}
        \caption{Time series of roll angle obtained by solving SDE (Eq.~(\ref{itopara})) }
        \label{fig:Timeseries_rollangle}
    \end{figure}

\begin{figure*}[tb]
        \centering
        \includegraphics[scale=0.8]{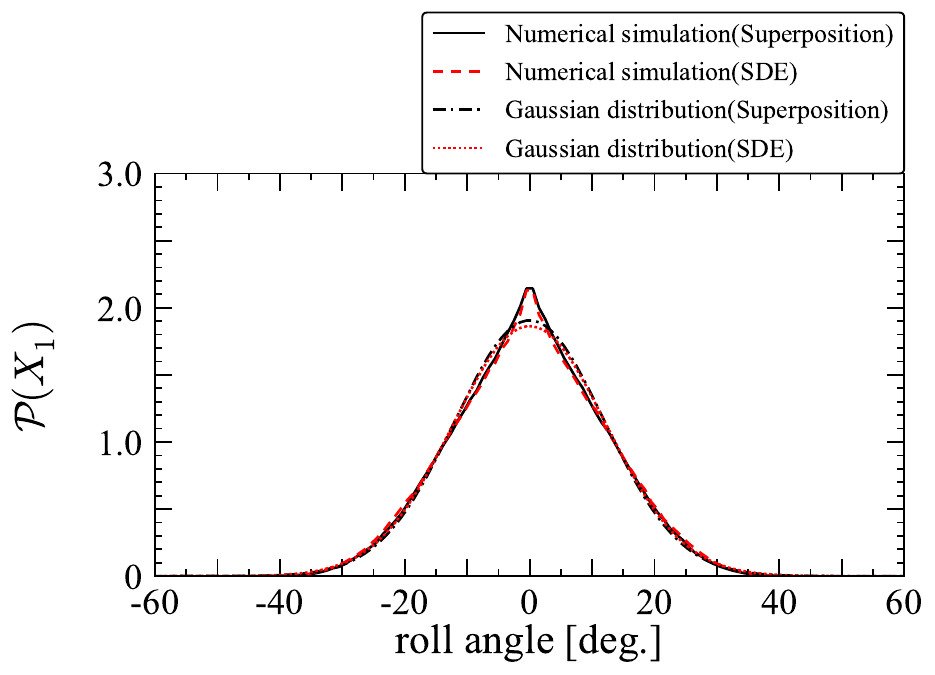}
        \includegraphics[scale=0.8]{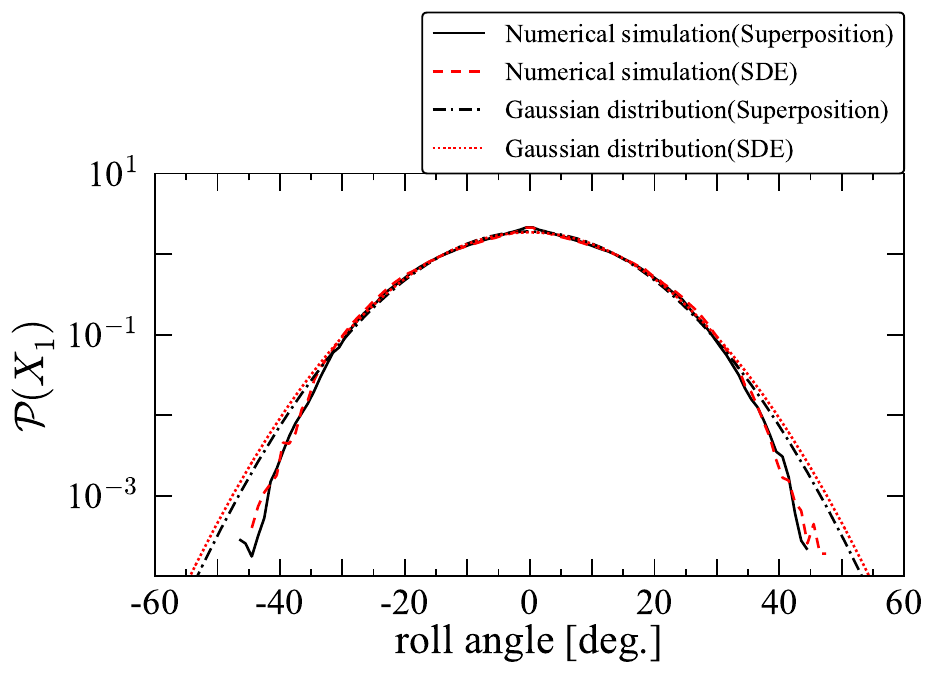}
        \caption{Comparison of PDF of $X_1$ among result obtained by solving SDE (Eq.~(\ref{itopara})) : Numerical simulation~(SDE), numerical simulation result by a principle of superposition : Numerical simulation~(Superposition), and Gaussian distributions by using moment values~(Table~\ref{tab:calculation_results_superposition_and_SDE}) }
        \label{fig:PDF_of_opt1}
    \end{figure*}

    As shown in Fig.~\ref{fig:Timeseries_rollangle}, it can be observed that the time series of roll angle, which is obtained by solving the SDE, shows a typical parametric rolling behaviour. As shown in Fig.~\ref{fig:PDF_of_opt1}, it can be observed that the PDF of roll angle obtained by solving the SDE has a satisfactory agreement with the PDF obtained by using the superposition principle. However, these PDFs do not agree with a Gaussian distribution.
    
\begin{figure*}[tb]
        \centering
        \includegraphics[scale=0.8]{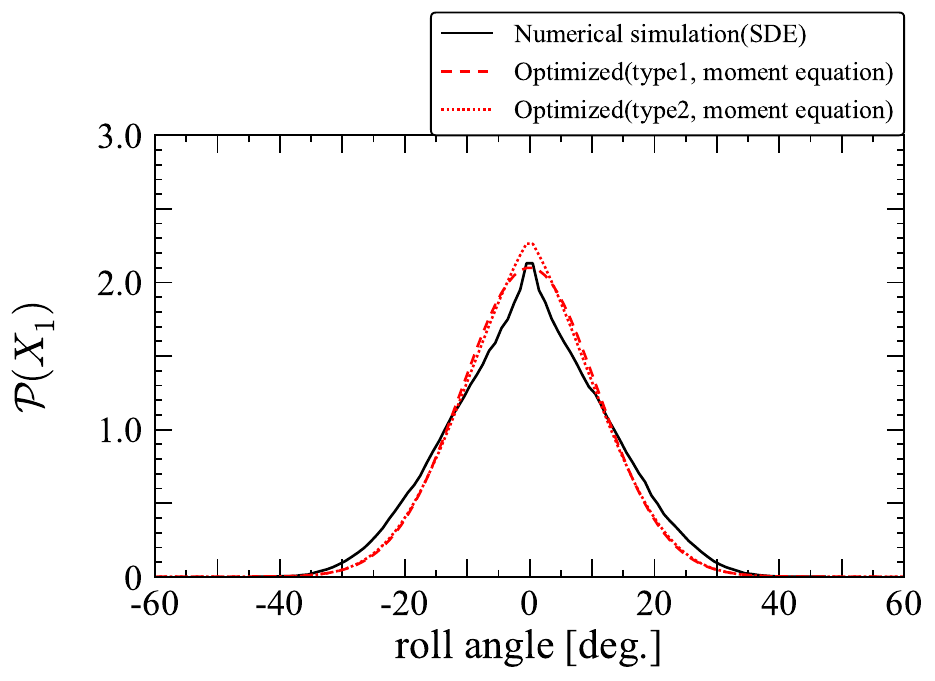}
        \includegraphics[scale=0.8]{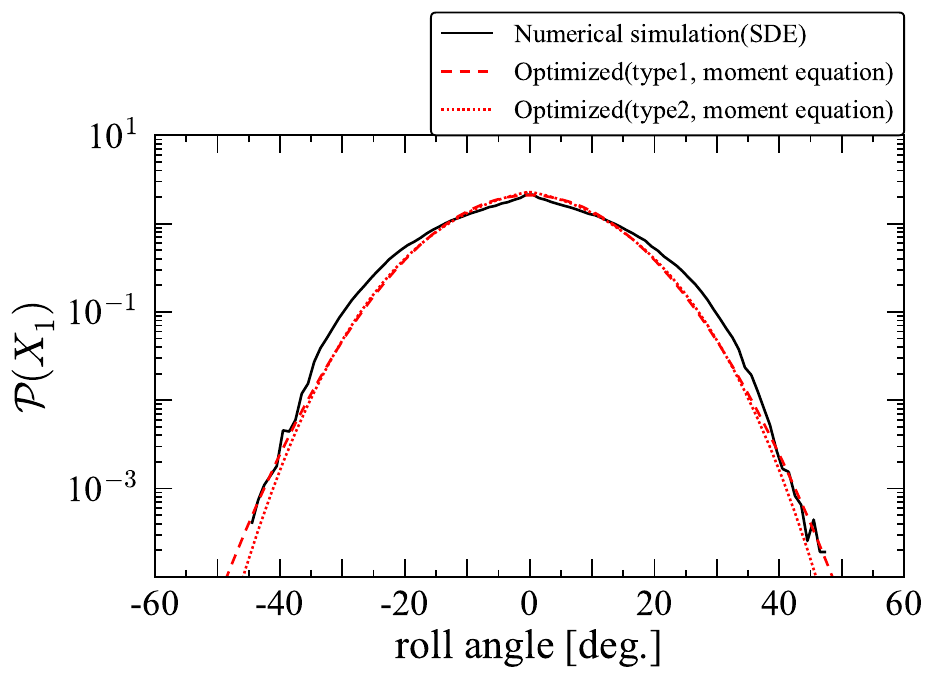}
        \caption{Comparison of PDF of $X_1$ among result obtained by solving SDE (Eq.~(\ref{itopara})) : Numerical simulation~(SDE), optimized result by using Eq.~(\ref{PDF_arbitrary_type1}) and moment values (Table~\ref{tab:calculation_results_moment_value_and_superposition}) obtained by solving moment equations (Eq.~(\ref{first_moment_para}) and (\ref{second_moment_para})) : optimized (type~1, moment equation), and optimized result by using Eq.~(\ref{PDF_arbitrary_type2}) and moment values (Table~\ref{tab:calculation_results_moment_value_and_superposition}) obtained by solving moment equations(Eq.~(\ref{first_moment_para}) and (\ref{second_moment_para})) : optimized (type~2, moment equation)}
        \label{fig:PDF_of_opt2}
    \end{figure*}

    Fig.~\ref{fig:PDF_of_opt2} shows the optimized result of PDFs of type~1 and type~2. Here, the moment values obtained by solving the moment equations are used. Applying the PDF of type~1, the obtained PDF has a similar function shape as a Gaussian distribution. This is because the second-order moment has a larger value than the other moments. We should focus on the point that the PDF of roll angle and the Laplace distribution resemble the function form. Therefore, in this study, Eq.~(\ref{PDF_arbitrary_type2}) is defined. As a result of using type~2, the optimized function has an improved agreement in shape with the numerical simulation result. 
    
\begin{figure*}[tb]
        \centering
        \includegraphics[scale=0.8]{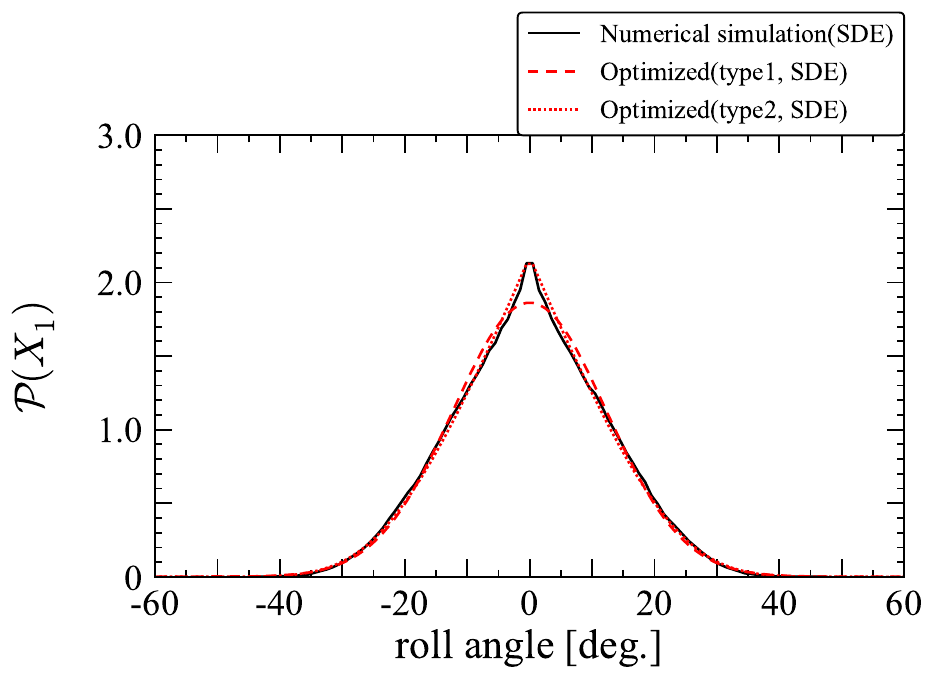}
        \includegraphics[scale=0.8]{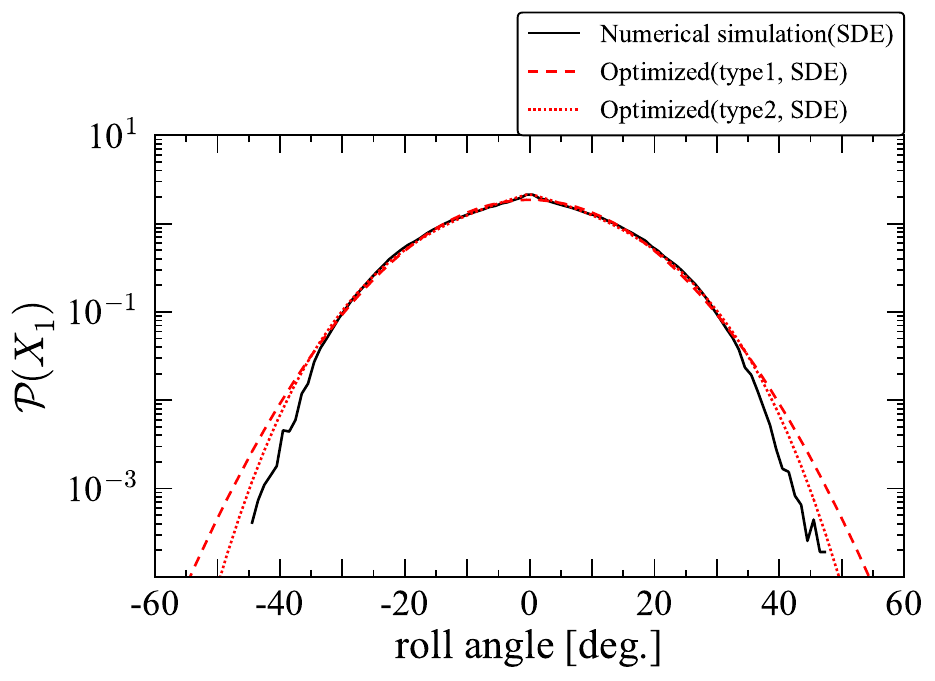}
        \caption{Comparison of PDF of $X_1$ among result obtained by solving SDE (Eq.~(\ref{itopara})) : Numerical simulation~(SDE), optimized result by using Eq.(\ref{PDF_arbitrary_type1}) and moment values (Table~\ref{tab:calculation_results_moment_value_and_superposition}) obtained by calculating SDE (Eq.~(\ref{itopara})) : optimized (type~1, SDE), and optimized result by using Eq.(\ref{PDF_arbitrary_type2}) and moment values(Table\ref{tab:calculation_results_moment_value_and_superposition}) obtained by calculating SDE (Eq.~(\ref{itopara})) : optimized (type~2, SDE)}
        \label{fig:PDF_of_opt3}
    \end{figure*}

    As we have seen, the moment values have a difference between the numerical result and the moment equation result. Thus, the moment values from the SDE results in Table~\ref{tab:calculation_results_superposition_and_SDE} are used. In Fig.~\ref{fig:PDF_of_opt3}, the optimized PDF of type~2 has a very good agreement with the SDE result. Therefore, it is clear that the optimized PDF, which has also a very good agreement with the numerical simulation result, can be obtained by using accurate moment values and the PDF of type~2.

    
\section{Concluding remarks}
It has been observed that the time series of the effective wave can be modeled sufficiently by using a linear filter. Furthermore, it turns out that the time series of GM variation can be modeled sufficiently by using the linear filter and the relationship between $\Delta$GM and wave height.

We have summarized a procedure in which the moment equations are obtained from a higher-order stochastic differential equation. Furthermore, we have suggested a simple way to calculate the steady state results for the moments. The appropriate moment values can be obtained easily for a higher-order SDE. Most important of all is the non-Gaussian PDF can be obtained by using the appropriate values for our proposed PDF and objective function.

In this study, the linear filter was used for the effective wave. We examine this filter whether it can be used for the methodology of GM variation and GZ variation. Further studies are needed in order to determine a set of probability density functions which have the shape of the objective probability densities.
    
\begin{acknowledgements}
This work was supported by a Grant-in-Aid for Scientific Research from the Japan Society for Promotion of Science (JSPS KAKENHI Grant Number 19H02360), as well as the collaborative research program / financial support from the Japan Society of Naval Architects and Ocean Engineers.
This study was supported by the Fundamental Research Developing Association for Shipbuilding and Offshore (REDAS), managed by the Shipbuilders’ Association of Japan from April 2020 to March 2023.
\end{acknowledgements}

\bibliographystyle{spphys}       
\bibliography{ref.bib}

\section*{Appendix 1}
The detailed formula about the relations between moments and cumulants can be shown in the electronic supplementary material (ESM).
    
\section*{Appendix 2}

In the case of using the second-order cumulant neglect closure, the stationary solutions of roll angle $X_{1}^{}$ and roll velocity $X_{2}^{}$ are obtained. Here, these solutions are oscillating, as shown in Fig.~\ref{fig:moment_value_in_unstable}. We consider the question whether the second-order cumulant method can sufficiently reflect the nonlinearity. Therefore, we examine the effect of considering higher-order cumulants. Thus, the third-order cumulant neglect closure is used as well. Thereby, one hundred twenty third-order moment equations, derived from Eq.~(\ref{moment_para}), are additionally needed. In other words, one hundred sixty-four moment equations are used in the corresponding numerical calculation. Furthermore, it is necessary that 4th and higher-order moments are represented by first, second, and third-order moments. These relations can be derived as shown in Appendix 1.

Fig.~\ref{fig:moment_value_in_unstable_third_order} shows the results when solving these moment equations numerically. Compared with the result of the second-order cumulant neglect closure in Fig.~\ref{fig:moment_value_in_unstable}, it is clear that the steady state solutions are not oscillating. Furthermore, the moment values in the third-order cumulant neglect closure are closer to the moment values obtained from the numerical result of the SDE than the second-order cumulant neglect closure results, as shown in Table~\ref{tab:calculation_results_moment_value_third_order}. We can therefore conclude that the result can reflect nonlinearity and be close to the actual value by using higher-order cumulants and moments.

    \renewcommand{\arraystretch}{1.5}
    \begin{table}[h]
      \begin{center}
        \caption{moment values obtained by solving moment equations (Eq.(\ref{first_moment_para}) and (\ref{second_moment_para}) )}
        \begin{tabular}{ll}
        \hline
         {Item} & {Value}
         \\
         \hline
          $\mathbb{E}\left[X_{1}^{}\right]$ & $0.00$
          \\
          $\mathbb{E}\left[X_{1}^{2}\right]$ & $3.72\times10_{}^{-2}$
          \\
          $\mathbb{E}\left[X_{2}^{}\right]$ & $0.00$
          \\
          $\mathbb{E}\left[X_{2}^{2}\right]$ & $2.42\times10_{}^{-3}$
          \\
        \hline
        \end{tabular}
        \label{tab:calculation_results_moment_value_third_order}
      \end{center}
    \end{table}
    \renewcommand{\arraystretch}{1.0}
    
     \begin{figure}[h]
        \centering
        \includegraphics[scale=1]{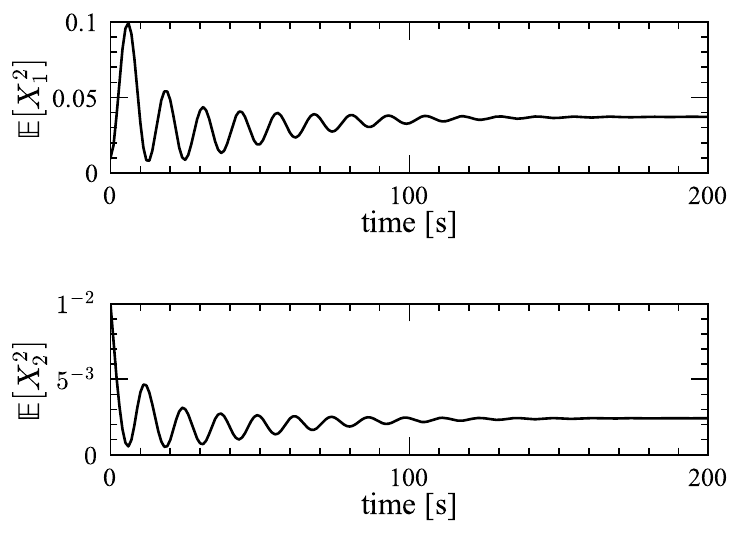}
        \caption{Second order moment of $X_{1}^{} \text{and} X_{2}^{}$ obtained by solving the ordinary differential equations of moment equation with third-order cumulant neglect closure method}
        \label{fig:moment_value_in_unstable_third_order}
    \end{figure}
    
\end{document}


\section*{Supplementary material}
In this study, the cumulant neglect closure method is used. We summarize the procedure for representing third and higher-order moments by first and second-order moments. 
Firstly, the relation expressions of cumulants are needed in order to consider the related expressions of the moments. Compared to the coefficients of the series expansion between moment generating function and cumulant generating function, the relating expressions between moments and cumulants can be obtained. Here, third and higher-order cumulants are set to zero.
The $k$th-order moment $m_{k_1^{}}$ of one variable is defined by
    \begin{equation}
        m_k^{} = \mathbb{E}[x_i^k] = \int_{-\,\infty}^{+\infty} x_i^k p(x_i^{})\,dx_i^{}\,\,.
    \end{equation}
$p(x_i^{})$ is a probability density function of one variable. Next, the moment~$m_{k_1^{},k_2^{},\cdots}^{}$( or $m_{k_1^{}k_2^{}\cdots}^{}$ ) of multivariable~$x_{i_1^{}},x_{i_2^{}}^{},\cdots$ is defined by
    \begin{equation}
        \begin{split}
        \displaystyle
            & m_{k_1^{},k_2^{},\cdots} = \mathbb{E}[x_{i_1^{}}^{k_1^{}}x_{i_2^{}}^{k_2^{}}\cdots] 
            \\
            \displaystyle
            & = \int_{-\,\infty}^{+\infty}\cdots\int_{-\,\infty}^{+\infty} x_{i_1^{}}^{k_1^{}} x_{i_2^{}}^{k_2^{}}\cdots p(x_{i_1^{}}^{},x_{i_2^{}}^{},\cdots)\,dx_{i_1^{}}^{}dx_{i_2^{}}^{}\cdots
        \end{split}
    \end{equation}
    $p(x_{i_1^{}},x_{i_2^{}},\cdots)$ is a multivariate probability density function. 
    
\,

\subsection*{Case1: Random variable is one.}\label{sec:appendix1case1}
In case that the random value is only one variable, the moment generating function $M(t)$ is represented as
    \begin{equation}
    \label{moment generating function}
        M(t)=\mathbb{E}[e^{tX}]=\sum_{i=0}^{\infty}\frac{\mathbb{E}[X^i]}{i!}t^i=1+t\mathbb{E}[X]+\frac{1}{2!}t^2\mathbb{E}[X^2]+\cdots\,\,.
    \end{equation}
    And the cumulant generating function is obtained by taking the logarithm of $M(t)$ as follows:
    \begin{equation}
    \label{cumulant generating function}
        \begin{split}
            \displaystyle 
            \log M(t) = \sum_{i=1}^{\infty}\frac{\kappa_{i}^{}}{i!}t^i
        \end{split}
        \,\,.
    \end{equation}
    Applying the Maclaurin expansion and the second-order cumulant neglect method~($\kappa_{i}^{}=0,\,i>2 $), 
    \begin{equation}
    \label{cumulant generating function 2}
        \begin{split}
        \displaystyle 
        & M(t)
        \\
        \displaystyle 
        = & \exp\left\{\sum_{i=1}^{2}\frac{\kappa_{i}^{}}{i!}t^i\right\}
        \\
        \displaystyle = & \exp(t\kappa_{1}^{})\exp\left(\frac{1}{2!}t^2\kappa_{2}^{}\right)
        \\
        \displaystyle = & \left\{ 1+t\kappa_{1}^{}+\frac{1}{2!}t^2\kappa_{1}^{2}+\frac{1}{3!}t^3\kappa_{1}^{3}+\cdots \right\}
        \\
        \displaystyle &
        \left\{ 1+\frac{1}{2!}t^2\kappa_{2}^{}+\frac{1}{2!}\left(\frac{1}{2!}t^2\kappa_{2}^{}\right)^2+\frac{1}{3!}\left(\frac{1}{2!}t^2\kappa_{2}^{}\right)^3+\cdots \right\}
        \\
        \displaystyle = & 1+t\kappa_{1}^{}+\frac{1}{2!}t^2(\kappa_{2}^{}+\kappa_{1}^{2})+\frac{1}{3!}t^3(\kappa_{1}^{3}+3\kappa_{1}^{}\kappa_{2}^{})+\cdots \,\,\,.
        \end{split}
    \end{equation}
    Comparing the coefficients between Eq.~(\ref{moment generating function}) and (\ref{cumulant generating function 2}), the first and second-order moments can be represented by means of the first and second-order cumulants.
    \begin{equation}
    \label{moment cumulant 1 and 2 order}
        \begin{split}
        \displaystyle m_{1}^{} = \kappa_{1}^{},\, m_{2}^{} = \kappa_{1}^{2}+\kappa_{2}^{}
        \end{split}
    \end{equation}
    Thereby, first and second-order cumulants can be represented by first and second-order moments.
    \begin{equation}
    \label{cumulant moment 1 and 2 order}
        \kappa_{1}^{} = m_{1}^{},\, \kappa_{2}^{} = m_{2}^{}-m_{1}^{2}
    \end{equation}
    Furthermore, substituting Eq.~(\ref{cumulant moment 1 and 2 order}) into the related expressions obtained by comparing the coefficients, the higher-order moments can be represented using the first and second-order moments. By this means Eq.~(\ref{moment cumulant higher order part1}) - (\ref{moment cumulant higher order part8}) are obtained.
    \begin{equation}
    \label{moment cumulant higher order part1}
        \begin{split}
            \displaystyle m_{3}^{} &= \kappa_{1}^{3} + 3\kappa_{1}^{}\kappa_{2}^{} 
            \\
            \displaystyle &= 3m_{1}^{}m_{2}^{} - 2m_{1}^{3}
        \end{split}
    \end{equation}
    \begin{equation}
    \label{moment cumulant higher order part2}
        \begin{split}
            \displaystyle m_{4}^{} &= \kappa_{1}^{4} + 6\kappa_{1}^{2}\kappa_{2}^{} + 3\kappa_{2}^{2} 
            \\
            \displaystyle &= 3m_{2}^{2} - 2m_{1}^{4}
        \end{split}
    \end{equation}
    \begin{equation}
    \label{moment cumulant higher order part3}
        \begin{split}
            \displaystyle m_{5}^{} &= \kappa_{1}^{5} + 10\kappa_{1}^{3}\kappa_{2}^{} + 15\kappa_{1}^{}\kappa_{2}^{2}
            \\
            \displaystyle &= 15m_{1}^{}m_{2}^{2} - 20m_{1}^{3}m_{2}^{} + 6m_{1}^{5}
        \end{split}
    \end{equation}    
    \begin{equation}
    \label{moment cumulant higher order part4}
        \begin{split}
            \displaystyle m_{6}^{} & = \kappa_{1}^{6} + 15\kappa_{1}^{4}\kappa_{2}^{} + 45\kappa_{1}^{2}\kappa_{2}^{2} + 15\kappa_{2}^{3}
            \\
            \displaystyle & = 15m_{2}^{3} - 30m_{1}^{4}m_{2}^{} + 16m_{1}^{6}
        \end{split}
    \end{equation}
    \begin{equation}
    \label{moment cumulant higher order part5}
        \begin{split}
            \displaystyle m_{7}^{} & = \kappa_{1}^{7} + 21\kappa_{1}^{5}\kappa_{2}^{} + 105\kappa_{1}^{3}\kappa_{2}^{2} + 105\kappa_{1}^{}\kappa_{2}^{3}
            \\
            \displaystyle & = 105m_{1}^{}m_{2}^{3} - 210m_{1}^{3}m_{2}^{2} + 126m_{1}^{5}m_{2}^{} -20m_{1}^{7}
        \end{split}
    \end{equation}
    \begin{equation}
    \label{moment cumulant higher order part6}
        \begin{split}
            \displaystyle m_{8}^{} & = \kappa_{1}^{8} + 28\kappa_{1}^{6}\kappa_{2}^{} + 210\kappa_{1}^{4}\kappa_{2}^{2} + 420\kappa_{1}^{2}\kappa_{2}^{3} + 105\kappa_{2}^{4}
            \\
            \displaystyle & = 105m_{2}^{4} - 420m_{1}^{4}m_{2}^{2} + 448m_{1}^{6}m_{2}^{} -132m_{1}^{8}
        \end{split}
    \end{equation}
    \begin{equation}
    \label{moment cumulant higher order part7}
        \begin{split}
            \displaystyle 
            m_{9}^{} & = \kappa_{1}^{9} + 36\kappa_{1}^{7}\kappa_{2}^{} + 378\kappa_{1}^{5}\kappa_{2}^{2} + 1260\kappa_{1}^{3}\kappa_{2}^{3} + 945\kappa_{1}^{}\kappa_{2}^{4}
            \\
            \displaystyle 
            & = 945m_{1}^{}m_{2}^{4} - 2520m_{1}^{3}m_{2}^{3} + 2268m_{1}^{5}m_{2}^{2} 
            \\
            \displaystyle 
            & -\, 720m_{1}^{7}m_{2}^{} + 28m_{1}^{9}
        \end{split}
    \end{equation}
    \begin{equation}
    \label{moment cumulant higher order part8}
        \begin{split}
            \displaystyle 
            m_{10}^{} & = \kappa_{1}^{10} + 45\kappa_{1}^{8}\kappa_{2}^{} + 630\kappa_{1}^{6}\kappa_{2}^{2} + 3150\kappa_{1}^{4}\kappa_{2}^{3} 
            \\
            \displaystyle 
            & +\, 4725\kappa_{1}^{2}\kappa_{2}^{4} + 945\kappa_{2}^{5}
            \\
            \displaystyle 
            & = 945m_{2}^{5} - 6300m_{1}^{4}m_{2}^{3} + 10080m_{1}^{6}m_{2}^{2}
            \\
            \displaystyle 
            & -\, 5940m_{1}^{8}m_{2}^{} + 1216m_{1}^{10}
        \end{split}
    \end{equation}

\subsection*{Case2:Random variables are two}\label{sec:appendix1case2}    
In case that the random values are two variables~$X_{1}^{},X_{2}^{}$, the moment generating function~$M(t_{1}^{},t_{2}^{})$ and the cumulant generating function are represented as in Eq.~(\ref{moment generating function 2 variable}) and (\ref{cumulant generating function 2 variable}).
    \begin{equation}
    \label{moment generating function 2 variable}
        \begin{split}
            \displaystyle
            M(t_{1}^{},t_{2}^{}) & = \mathbb{E}[e^{t_{1}^{}X_{1}^{}}e^{t_{2}^{}X_{2}^{}}]
            = 1 + \sum_{i=1}^{\infty}\frac{\mathbb{E}[X_{1}^{i}]}{i!}t_{1}^{i} 
            \\
            \displaystyle 
            & + \sum_{j=1}^{\infty}\frac{\mathbb{E}[X_{2}^{j}]}{j!}t_{2}^{j} + \sum_{i=1}^{\infty}\sum_{j=1}^{\infty}\frac{\mathbb{E}[X_{1}^{i}X_{2}^{j}]}{i!j!}t_{1}^{i}t_{2}^{j}
        \end{split}
    \end{equation}
    \begin{equation}
    \label{cumulant generating function 2 variable}
        \begin{split}
            \displaystyle
            \log M(t_{1}^{},t_{2}^{})
            = \sum_{i=1}^{\infty}\frac{\kappa_{i\,0}^{}}{i!}t_{1}^{i} + \sum_{j=1}^{\infty}\frac{\kappa_{0\,j}^{}}{j!}t_{2}^{j}
            + \sum_{i=1}^{\infty}\sum_{j=1}^{\infty}\frac{\kappa_{i\,j}^{}}{i!j!}t_{1}^{i}t_{2}^{j}
        \end{split}
    \end{equation}
    Applying the Maclaurin expansion and the second-order cumulant neglect method~($\kappa_{ij}^{}=0,\,i+j>2$), 
    \begin{equation}
    \label{cumulant generating function 2 variable 2}
        \begin{split}
            \displaystyle
            M(t_{1}^{},t_{2}^{})
            = & \exp\left( t_{1}^{}\kappa_{10}^{} \right) \exp\left( \frac{t_{1}^{2}}{2!}\kappa_{20}^{} \right) \exp\left( t_{2}^{}\kappa_{01}^{} \right)
            \\
            \displaystyle &
            \exp\left( \frac{t_{2}^{2}}{2!}\kappa_{02}^{} \right) \exp\left( t_{1}^{}t_{2}^{}\kappa_{11}^{} \right)
            \\
            \displaystyle
            = &\left\{ 1+t_{1}^{}\kappa_{10}^{}+\frac{1}{2!}t_{1}^{2}\kappa_{20}^{2}+\cdots \right\}
            \\
            \displaystyle &
            \left\{ 1+\frac{1}{2!}t_{1}^{2}\kappa_{20}^{}+\frac{1}{2!}\left(\frac{1}{2!}t_{1}^{2}\kappa_{20}^{}\right)^2+\cdots \right\}
            \\
            \displaystyle &
            \left\{ 1+t_{2}^{}\kappa_{01}^{}+\frac{1}{2!}t_{2}^{2}\kappa_{01}^{2}+\cdots \right\}
            \\
            \displaystyle &
            \left\{ 1+\frac{1}{2!}t_{2}^{2}\kappa_{02}^{}+\frac{1}{2!}\left(\frac{1}{2!}t_{2}^{2}\kappa_{02}^{}\right)^2+\cdots \right\}
            \\
            \displaystyle &
            \left\{ 1+t_{1}^{}t_{2}^{}\kappa_{11}^{}+\frac{1}{2!}t_{1}^{2}t_{2}^{2}\kappa_{11}^{2}+\cdots \right\}\,\,.
        \end{split}
    \end{equation}
    Comparing the coefficients between Eq.~(\ref{moment generating function 2 variable}) and (\ref{cumulant generating function 2 variable 2}), the first and second-order moments can be represented using the first and second-order cumulants.
    \begin{equation}
    \label{moment cumulant 1 and 2 order 2 variable}
        \begin{split}
            \displaystyle 
            & m_{10}^{} = \kappa_{10}^{},\, m_{01}^{} = \kappa_{01}^{},\, m_{11}^{} = \kappa_{01}^{}\kappa_{10}^{}+\kappa_{11}^{}
            \\
            \displaystyle & m_{20}^{} = \kappa_{10}^{2} + \kappa_{20}^{},\, m_{02}^{} = \kappa_{01}^{2} + \kappa_{02}^{}
        \end{split}
    \end{equation}
    Thereby, the first and second-order cumulants can be represented using the first and second-order moments.
    \begin{equation}
    \label{cumulant moment 1 and 2 order 2 variable}
        \begin{split}
            \displaystyle 
            & \kappa_{10}^{} = m_{10}^{},\, \kappa_{01}^{} = m_{01}^{},\, \kappa_{11}^{} = m_{11}^{} - m_{10}^{}m_{01}^{}
            \\
            \displaystyle 
            & \kappa_{20}^{} = m_{20}^{} - m_{10}^{2},\, \kappa_{02}^{} = m_{02}^{} - m_{01}^{2}
        \end{split}
    \end{equation}
    Furthermore, substituting Eq.~(\ref{cumulant moment 1 and 2 order 2 variable}) into the related expressions obtained by coefficient comparison, the higher-order moments can be represented using the first and second-order moments. By this means Eq.~(\ref{moment cumulant higher order 2 variable part1}) - (\ref{moment cumulant higher order 2 variable part11}) are obtained.
    \begin{equation}
    \label{moment cumulant higher order 2 variable part1}
        \begin{split}
            \displaystyle m_{1,2}^{} & = \kappa_{01}^{2}\kappa_{10}^{} + \kappa_{02}^{}\kappa_{10}^{} + 2\kappa_{01}^{}\kappa_{11}^{}
            \\
            \displaystyle  & = m_{10}^{}m_{02}^{} - 2m_{01}^{2}m_{10}^{} + 2m_{01}^{}m_{11}^{}
        \end{split}
    \end{equation}
    \begin{equation}
    \label{moment cumulant higher order 2 variable part2}
        \begin{split}
            \displaystyle m_{1,3}^{} & = \kappa_{01}^{3}\kappa_{10}^{} + 3\kappa_{02}^{}\kappa_{10}^{}\kappa_{01}^{} + 3\kappa_{01}^{2}\kappa_{11}^{} + 3\kappa_{02}^{}\kappa_{11}^{}
            \\
            \displaystyle & = 3m_{11}^{}m_{02}^{} - 2m_{01}^{3}m_{10}^{}
        \end{split}
    \end{equation}    
    \begin{equation}
    \label{moment cumulant higher order 2 variable part3}
        \begin{split}
            \displaystyle
            m_{1,4}^{} & = \kappa_{10}^{} \kappa_{01}^4 + 4 \kappa_{11}^{} \kappa_{01}^3 + 6 \kappa_{02}^{} \kappa_{10}^{} \kappa_{01}^2
            \\
            \displaystyle &
            +\, 12 \kappa_{02}^{} \kappa_{11}^{}\kappa_{01}^{}+3 \kappa_{02}^2 \kappa_{10}^{}
            \\
            \displaystyle & = 6m_{01}^{4}m_{10}^{} - 12m_{01}^{2}m_{10}^{}m_{02}^{} + 3m_{02}^{2}m_{10}^{}
            \\
            \displaystyle &
            -\, 8m_{01}^{3}m_{11}^{} + 12m_{01}^{}m_{11}^{}m_{02}^{}
        \end{split}
    \end{equation}    
    \begin{equation}
    \label{moment cumulant higher order 2 variable part4}
        \begin{split}
            \displaystyle m_{1,5}^{} & = \kappa_{01}^{5}\kappa_{10}^{} + 10\kappa_{02}^{}\kappa_{10}^{}\kappa_{01}^{2} + 15\kappa_{01}^{}\kappa_{02}^{2}\kappa_{10}^{} 
            \\
            \displaystyle &
            +\, 5\kappa_{01}^{4}\kappa_{11}^{} + 30\kappa_{02}^{}\kappa_{11}^{}\kappa_{01}^{2} + 15\kappa_{11}^{}\kappa_{02}^{2}
            \\
            \displaystyle & = 16m_{01}^{5}m_{10}^{} - 20m_{01}^{3}m_{10}^{}m_{02}^{} - 10m_{01}^{4}m_{11}^{}
            \\
            \displaystyle &
            +\, 15m_{02}^{2}m_{11}^{}
        \end{split}
    \end{equation}    
    \begin{equation}
    \label{moment cumulant higher order 2 variable part5}
        \begin{split}
            \displaystyle m_{1,6}^{} & = \kappa_{01}^{6}\kappa_{10}^{} + 15\kappa_{01}^{4}\kappa_{10}^{}\kappa_{02}^{} + 45\kappa_{01}^{2}\kappa_{02}^{2}\kappa_{10}^{}
            \\
            \displaystyle &
            +\, 15\kappa_{02}^{3}\kappa_{10}^{} + 6\kappa_{01}^{5}\kappa_{11}^{} + 60\kappa_{01}^{3}\kappa_{02}^{}\kappa_{11}^{}
            \\
            \displaystyle &
            +\, 90\kappa_{01}^{}\kappa_{02}^{2}\kappa_{11}^{}
            \\
            \displaystyle & = - 20m_{01}^{6}m_{10}^{} + 90m_{01}^{4}m_{10}^{}m_{02}^{} - 90m_{01}^{2}m_{02}^{2}m_{10}^{}
            \\
            \displaystyle &
            +\, 15m_{02}^{3}m_{10}^{} + 36m_{01}^{5}m_{11}^{} - 120m_{01}^{3}m_{11}^{}m_{02}^{}
            \\
            \displaystyle &
            +\, 90m_{01}^{}m_{11}^{}m_{02}^{2}
        \end{split}
    \end{equation}    
    \begin{equation}
    \label{moment cumulant higher order 2 variable part6}
        \begin{split}
            \displaystyle m_{1,7}^{} & = \kappa_{01}^{7}\kappa_{10}^{} + 21\kappa_{01}^{5}\kappa_{10}^{}\kappa_{02}^{} + 105\kappa_{01}^{3}\kappa_{02}^{2}\kappa_{10}^{} 
            \\
            \displaystyle &
            +\, 7\kappa_{01}^{6}\kappa_{11}^{} + 105\kappa_{02}^{3}\kappa_{10}^{}\kappa_{01}^{} + 105\kappa_{02}^{3}\kappa_{11}^{}
            \\
            \displaystyle &
            +\, 105\kappa_{01}^{4}\kappa_{02}^{}\kappa_{11}^{} + 315\kappa_{01}^{2}\kappa_{02}^{2}\kappa_{11}^{}
            \\
            \displaystyle & = - 132m_{01}^{7}m_{10}^{} + 336m_{01}^{5}m_{10}^{}m_{02}^{}
            \\
            \displaystyle &
            -\, 210m_{01}^{3}m_{02}^{2}m_{10}^{} + 112m_{01}^{6}m_{11}^{}
            \\
            \displaystyle &
            -\, 210m_{01}^{4}m_{11}^{}m_{02}^{} + 105m_{11}^{}m_{02}^{3}
        \end{split}
    \end{equation}    
    \begin{equation}
    \label{moment cumulant higher order 2 variable part7}
        \begin{split}
            \displaystyle m_{1,8}^{} & = \kappa_{01}^{8}\kappa_{10}^{} + 28\kappa_{01}^{6}\kappa_{10}^{}\kappa_{02}^{} + 210\kappa_{01}^{4}\kappa_{02}^{2}\kappa_{10}^{} 
            \\
            \displaystyle &
            +\, 420\kappa_{02}^{3}\kappa_{10}^{}\kappa_{01}^{2} + 105\kappa_{02}^{4}\kappa_{10}^{} + 168\kappa_{01}^{5}\kappa_{02}^{}\kappa_{11}^{}
            \\
            \displaystyle &
            + 840\kappa_{01}^{3}\kappa_{02}^{2}\kappa_{11}^{} + 840\kappa_{01}^{}\kappa_{02}^{3}\kappa_{11}^{} + 8\kappa_{01}^{7}\kappa_{11}^{}
            \\
            \displaystyle & = 28m_{01}^{8}m_{10}^{} - 560m_{01}^{6}m_{10}^{}m_{02}^{} + 1260m_{01}^{4}m_{02}^{2}m_{10}^{}
            \\
            \displaystyle &
            -\, 840m_{01}^{2}m_{10}^{}m_{02}^{3} + 105m_{10}^{}m_{02}^{4} - 160m_{11}^{}m_{01}^{7} 
            \\
            \displaystyle &
            +\, 1008m_{01}^{5}m_{11}^{}m_{02}^{} - 1680m_{01}^{3}m_{02}^{2}m_{11}^{} 
            \\
            \displaystyle &
            +\, 840m_{01}^{}m_{02}^{3}m_{11}^{}
        \end{split}
    \end{equation}    
    \begin{equation}
    \label{moment cumulant higher order 2 variable part8}
        \begin{split}
        \displaystyle
            m_{1,9}^{} & = \kappa_{01}^{9}\kappa_{10}^{} + 36\kappa_{01}^{7}\kappa_{10}^{}\kappa_{02}^{} + 378\kappa_{01}^{5}\kappa_{10}^{}\kappa_{02}^{2}
            \\
            \displaystyle &
            +\, 9\kappa_{11}^{}\kappa_{01}^{8} + 945\kappa_{01}^{}\kappa_{02}^{4}\kappa_{10}^{} + 252\kappa_{01}^{6}\kappa_{02}^{}\kappa_{11}^{}
            \\
            \displaystyle &
            +\, 1260\kappa_{01}^{3}\kappa_{10}^{}\kappa_{02}^{3} + 1890\kappa_{01}^{4}\kappa_{02}^{2}\kappa_{11}^{}
            \\
            \displaystyle &
            +\, 3780\kappa_{01}^{2}\kappa_{02}^{3}\kappa_{11}^{} + 945\kappa_{02}^{4}\kappa_{11}^{}
            \\
            \displaystyle &
            = 1216m_{01}^{9}m_{10}^{} - 4752m_{01}^{7}m_{10}^{}m_{02}^{} 
            \\
            \displaystyle &
            +\, 6048m_{01}^{5}m_{02}^{2}m_{10}^{} - 2520m_{01}^{3}m_{10}^{}m_{02}^{3} 
            \\
            \displaystyle &
            -\, 1188m_{01}^{8}m_{11}^{} + 945m_{02}^{4}m_{11}^{}
            \\
            \displaystyle &
            - 3780m_{01}^{4}m_{11}^{}m_{02}^{2} + 4032m_{11}^{}m_{01}^{6}m_{02}^{}
        \end{split}
    \end{equation}    
    \begin{equation}
    \label{moment cumulant higher order 2 variable part9}
        \begin{split}
            \displaystyle
            m_{1,10}^{} & = \kappa_{10}^{} \kappa_{01}^{10} + 10 \kappa_{11}^{} \kappa_{01}^9 + 45 \kappa_{02}^{} \kappa_{10}^{} \kappa_{01}^8 
            \\
            \displaystyle &
            +\, 360 \kappa_{02}^{} \kappa_{11}^{} \kappa_{01}^7 + 630 \kappa_{02}^2 \kappa_{10}^{} \kappa_{01}^6
            \\
            \displaystyle &
            +\, 3780 \kappa_{02}^2 \kappa_{11}^{} \kappa_{01}^5 + 3150 \kappa_{02}^3 \kappa_{10}^{} \kappa_{01}^4
            \\
            \displaystyle &
            +\, 12600 \kappa_{02}^3 \kappa_{11}^{} \kappa_{01}^3 + 4725 \kappa_{02}^4 \kappa_{10}^{} \kappa_{01}^2 
            \\
            \displaystyle &
            +\, 9450 \kappa_{02}^4 \kappa_{11}^{} \kappa_{01}^{} + 945 \kappa_{02}^5 \kappa_{10}^{}
            \\
            \displaystyle &
            = 936 m_{10}^{} m_{01}^{10} + 280 m_{11}^{} m_{01}^9 + 1260 m_{02}^{} m_{10}^{}m_{01}^8
            \\
            \displaystyle &
            -\, 7200 m_{02}^{} m_{11}^{} m_{01}^7 - 12600 m_{02}^2 m_{10}^{} m_{01}^6 
            \\
            \displaystyle &
            +\, 22680 m_{02}^2 m_{11}^{} m_{01}^5 + 18900 m_{02}^3 m_{10}^{} m_{01}^4
            \\
            \displaystyle &
            -\, 25200 m_{02}^3 m_{11}^{} m_{01}^3 - 9450 m_{02}^4 m_{10}^{} m_{01}^2
            \\
            \displaystyle &
            +\, 9450 m_{02}^4 m_{11}^{} m_{01}^{} + 945 m_{02}^5 m_{10}^{}
        \end{split}
    \end{equation}
    \begin{equation}
    \label{moment cumulant higher order 2 variable part10}
        \begin{split}
            \displaystyle 
            m_{1,11}^{} & = \kappa_{10}^{} \kappa_{01}^{11} + 11 \kappa_{11}^{} \kappa_{01}^{10} + 55 \kappa_{02}^{} \kappa_{10}^{} \kappa_{01}^9 
            \\
            \displaystyle &
            +\, 495 \kappa_{02}^{} \kappa_{11}^{} \kappa_{01}^8 + 990 \kappa_{02}^2 \kappa_{10}^{} \kappa_{01}^7
            \\
            \displaystyle &
            +\, 6930 \kappa_{02}^2 \kappa_{11}^{} \kappa_{01}^6 + 6930 \kappa_{02}^3 \kappa_{10}^{} \kappa_{01}^5
            \\
            \displaystyle &
            +\, 34650 \kappa_{02}^3 \kappa_{11}^{} \kappa_{01}^4 + 17325 \kappa_{02}^4 \kappa_{10}^{} \kappa_{01}^3 
            \\
            \displaystyle &
            +\, 51975 \kappa_{02}^4 \kappa_{11}^{} \kappa_{01}^2 + 10395 \kappa_{02}^5 \kappa_{10}^{} \kappa_{01}^{}
            \\
            \displaystyle &
            +\, 10395 \kappa_{02}^5 \kappa_{11}^{}
            \\
            \displaystyle &
            = -12440 m_{10}^{} m_{01}^{11} + 13376 m_{11}^{} m_{01}^{10} 
            \\
            \displaystyle &
            +\, 66880 m_{02}^{} m_{10}^{}m_{01}^9 - 65340 m_{02}^{} m_{11}^{} m_{01}^8
            \\
            \displaystyle &
            -\, 130680 m_{02}^2 m_{10}^{}m_{01}^7 + 110880 m_{02}^2 m_{11}^{} m_{01}^6 
            \\
            \displaystyle &
            +\, 110880 m_{02}^3 m_{10}^{}m_{01}^5 - 69300 m_{02}^3 m_{11}^{} m_{01}^4
            \\
            \displaystyle &
            - 34650 m_{02}^4 m_{10}^{}m_{01}^3 + 10395 m_{02}^5 m_{11}^{}
        \end{split}
    \end{equation} 
    \begin{equation}
    \label{moment cumulant higher order 2 variable part11}
        \begin{split}
        \displaystyle
            m_{1,12}^{} & = \kappa_{10}^{} \kappa_{01}^{12} + 12 \kappa_{11}^{} \kappa_{01}^{11} + 66 \kappa_{02}^{} \kappa_{10}^{}\kappa_{01}^{10}
            \\
            \displaystyle &
            +\, 660 \kappa_{02}^{} \kappa_{11}^{} \kappa_{01}^9 + 1485 \kappa_{02}^2 \kappa_{10}^{}\kappa_{01}^8 
            \\
            \displaystyle &
            +\, 11880 \kappa_{02}^2 \kappa_{11}^{} \kappa_{01}^7 + 13860 \kappa_{02}^3 \kappa_{10}^{}\kappa_{01}^6
            \\
            \displaystyle &
            +\, 83160 \kappa_{02}^3 \kappa_{11}^{} \kappa_{01}^5 + 51975 \kappa_{02}^4 \kappa_{10}^{}\kappa_{01}^4 
            \\
            \displaystyle &
            +\, 207900 \kappa_{02}^4 \kappa_{11}^{} \kappa_{01}^3 + 62370 \kappa_{02}^5 \kappa_{10}^{}\kappa_{01}^2 
            \\
            \displaystyle &
            +\, 124740 \kappa_{02}^5 \kappa_{11}^{} \kappa_{01}^{} + 10395 \kappa_{02}^6 \kappa_{10}^{}
            \\
            \displaystyle &
            = -23672 m_{10}^{} m_{01}^{12} + 11232 m_{11}^{} m_{01}^{11} 
            \\
            \displaystyle &
            +\, 61776 m_{02}^{} m_{10}^{}m_{01}^{10} + 18480 m_{02}^{} m_{11}^{} m_{01}^9 
            \\
            \displaystyle &
            +\, 41580 m_{02}^2 m_{10}^{}m_{01}^8 - 237600 m_{02}^2 m_{11}^{} m_{01}^7 
            \\
            \displaystyle &
            -\, 277200 m_{02}^3 m_{10}^{}m_{01}^6 + 498960 m_{02}^3 m_{11}^{} m_{01}^5 
            \\
            \displaystyle &
            +\, 311850 m_{02}^4 m_{10}^{}m_{01}^4 - 415800 m_{02}^4 m_{11}^{}m_{01}^3
            \\
            \displaystyle &
            -\, 124740 m_{02}^5 m_{10}^{}m_{01}^2 + 124740 m_{02}^5 m_{11}^{} m_{01}^{}
            \\
            \displaystyle &
            +\, 10395 m_{02}^6 m_{10}^{}
        \end{split}
    \end{equation}

\subsection*{Case2:Random variables are three}\label{sec:appendix1case3} 
    In case of three random variables~$X_{1}^{},X_{2}^{},X_{3}^{}$, the moment generating function~$M(t_{1}^{},t_{2}^{},t_{3}^{})$ and the cumulant generating function are represented as in Eqs.~(\ref{moment generating function 3 variable}) and (\ref{cumulant generating function 3 variable}).
    \begin{equation}
    \label{moment generating function 3 variable}
        \begin{split}
            \displaystyle & M(t_{1}^{},t_{2}^{},t_{3}^{})
            \\
            \displaystyle 
            = & \mathbb{E}[e^{t_{1}^{}X_{1}^{}}e^{t_{2}^{}X_{2}^{}}e^{t_{3}^{}X_{3}^{}}]
            \\
            \displaystyle
            = & 1 + \sum_{i=1}^{\infty}\frac{\mathbb{E}[X_{1}^{i}]}{i!}t_{1}^{i}
            + \sum_{j=1}^{\infty}\frac{\mathbb{E}[X_{2}^{j}]}{j!}t_{2}^{j}
            + \sum_{k=1}^{\infty}\frac{\mathbb{E}[X_{3}^{k}]}{k!}t_{3}^{k}
            \\
            \displaystyle
            + & \sum_{i=1}^{\infty}\sum_{j=1}^{\infty}\frac{\mathbb{E}[X_{1}^{i}X_{2}^{j}]}{i!j!}t_{1}^{i}t_{2}^{j}
            + \sum_{i=1}^{\infty}\sum_{k=1}^{\infty}\frac{\mathbb{E}[X_{1}^{i}X_{3}^{k}]}{i!k!}t_{1}^{i}t_{3}^{k}
            \\
            \displaystyle
            + & \sum_{j=1}^{\infty}\sum_{k=1}^{\infty}\frac{\mathbb{E}[X_{2}^{j}X_{3}^{k}]}{j!k!}t_{2}^{j}t_{3}^{k}
            + \sum_{i=1}^{\infty}\sum_{j=1}^{\infty}\sum_{k=1}^{\infty}\frac{\mathbb{E}[X_{1}^{i}X_{2}^{j}X_{3}^{k}]}{i!j!k!}t_{1}^{i}t_{2}^{j}t_{3}^{k}
        \end{split}
    \end{equation}
    \begin{equation}
    \label{cumulant generating function 3 variable}
        \begin{split}
            \displaystyle
            &\log M(t_{1}^{},t_{2}^{},t_{3}^{})
            \\
            \displaystyle
            = & \sum_{i=1}^{\infty}\frac{\kappa_{i\,0\,0}^{}}{i!}t_{1}^{i} + \sum_{j=1}^{\infty}\frac{\kappa_{0\,j\,0}^{}}{j!}t_{2}^{j} + \sum_{k=1}^{\infty}\frac{\kappa_{0\,0\,k}^{}}{k!}t_{3}^{k}
            \\
            \displaystyle
            + & \sum_{i=1}^{\infty}\sum_{j=1}^{\infty}\frac{\kappa_{i\,j\,0}^{}}{i!j!}t_{1}^{i}t_{2}^{j}
            + \sum_{i=1}^{\infty}\sum_{k=1}^{\infty}\frac{\kappa_{i\,0\,k}^{}}{i!k!}t_{1}^{i}t_{3}^{k}
            + \sum_{j=1}^{\infty}\sum_{k=1}^{\infty}\frac{\kappa_{0\,j\,k}^{}}{j!k!}t_{2}^{j}t_{3}^{k}
            \\
            \displaystyle
            + & \sum_{i=1}^{\infty}\sum_{j=1}^{\infty}\sum_{k=1}^{\infty}\frac{\kappa_{i\,j\,k}^{}}{i!j!k!}t_{1}^{i}t_{2}^{j}t_{3}^{k}
        \end{split}
    \end{equation}    
    Applying the Maclaurin expansion and the second-order cumulant neglect method~($\kappa_{ijk}^{}=0,\,i+j+k>2$), 
    \begin{equation}
    \label{cumulant generating function 3 variable 2}
        \begin{split}
            \displaystyle
            &M(t_{1}^{},t_{2}^{},t_{3}^{})
            \\
            \displaystyle
            = & \exp\left( t_{1}^{}\kappa_{100}^{} \right) \exp\left( \frac{t_{1}^{2}}{2!}\kappa_{200}^{} \right) \exp\left( t_{2}^{}\kappa_{010}^{} \right)
            \\
            \displaystyle &
            \exp\left( \frac{t_{2}^{2}}{2!}\kappa_{020}^{} \right) \exp\left( t_{3}^{}\kappa_{001}^{} \right) \exp\left( \frac{t_{3}^{2}}{2!}\kappa_{002}^{} \right)
            \\
            \displaystyle &
            \exp\left( t_{1}^{}t_{2}^{}\kappa_{110}^{} \right) \exp\left( t_{1}^{}t_{3}^{}\kappa_{101}^{} \right) \exp\left( t_{2}^{}t_{3}^{}\kappa_{011}^{} \right)
            \\
            \displaystyle 
            = & \left\{ 1+t_{1}^{}\kappa_{100}^{} +\frac{1}{2!}t_{1}^{2}\kappa_{100}^{2}+\cdots \right\}
            \\
            \displaystyle &
            \left\{ 1+\frac{1}{2!}t_{1}^{2}\kappa_{200}^{}+\frac{1}{2!}\left(\frac{1}{2!}t_{1}^{2}\kappa_{200}^{}\right)^2+\cdots \right\}
            \\
            \displaystyle &
            \left\{ 1+t_{2}^{}\kappa_{010}^{}+\frac{1}{2!}t_{2}^{2}\kappa_{010}^{2}+\cdots \right\}
            \\
            \displaystyle &
            \left\{ 1+\frac{1}{2!}t_{2}^{2}\kappa_{020}^{}+\frac{1}{2!}\left(\frac{1}{2!}t_{2}^{2}\kappa_{020}^{}\right)^2+\cdots \right\}
            \\
            \displaystyle &
            \left\{ 1+t_{3}^{}\kappa_{001}^{}+\frac{1}{2!}t_{3}^{2}\kappa_{001}^{2}+\cdots \right\}
            \\
            \displaystyle &
            \left\{ 1+\frac{1}{2!}t_{3}^{2}\kappa_{002}^{}+\frac{1}{2!}\left(\frac{1}{2!}t_{3}^{2}\kappa_{002}^{}\right)^2+\cdots \right\}
            \\
            \displaystyle &
            \left\{ 1+t_{1}^{}t_{2}^{}\kappa_{110}^{}+\frac{1}{2!}t_{1}^{2}t_{2}^{2}\kappa_{110}^{2}+\cdots \right\}
            \\
            \displaystyle &
            \left\{ 1+t_{1}^{}t_{3}^{}\kappa_{101}^{}+\frac{1}{2!}t_{1}^{2}t_{3}^{2}\kappa_{101}^{2}+\cdots \right\}
            \\
            \displaystyle &
            \left\{ 1+t_{2}^{}t_{3}^{}\kappa_{011}^{}+\frac{1}{2!}t_{2}^{2}t_{3}^{2}\kappa_{011}^{2}+\cdots \right\}\,\,.
        \end{split}
    \end{equation}
    Comparing the coefficients between Eq.~(\ref{moment generating function 3 variable}) and (\ref{cumulant generating function 3 variable 2}), the first and second-order moments can be represented using the first and second-order cumulants.
    \begin{equation}
    \label{moment cumulant 1 and 2 order 3 variable}
        \begin{split}
            \displaystyle &
            m_{100}^{} = \kappa_{100}^{},\, m_{010}^{} = \kappa_{010}^{},\, m_{001}^{} = \kappa_{001}^{}
            \\
            \displaystyle &
            m_{110}^{} = \kappa_{100}^{}\kappa_{010}^{} + \kappa_{110}^{},\, m_{101}^{} = \kappa_{100}^{}\kappa_{001}^{} + \kappa_{101}^{}
            \\
            \displaystyle &
            m_{011}^{} = \kappa_{010}^{}\kappa_{001}^{} + \kappa_{011}^{},\, m_{200}^{} = \kappa_{100}^{2} + \kappa_{200}^{}
            \\
            \displaystyle &
            m_{020}^{} = \kappa_{010}^{2} + \kappa_{020}^{},\, m_{002}^{} = \kappa_{001}^{2} + \kappa_{002}^{}
        \end{split}
    \end{equation}
    Thereby, the first and second-order cumulants can be represented using the first and second-order moments.
    \begin{equation}
    \label{cumulant moment 1 and 2 order 3 variable}
        \begin{split}
            \displaystyle &
            \kappa_{100}^{} = m_{100}^{},\, \kappa_{010}^{} = m_{010}^{},\, \kappa_{001}^{} = m_{001}^{}
            \\
            \displaystyle &
            \kappa_{110}^{} = m_{110}^{} - m_{100}^{}m_{010}^{},\, \kappa_{101}^{} = m_{101}^{} - m_{100}^{}m_{001}^{}
            \\
            \displaystyle &
            \kappa_{011}^{} = m_{011}^{} - m_{010}^{}m_{001}^{},\, \kappa_{200}^{} = m_{200}^{} - m_{100}^{2}
            \\
            \displaystyle &
            \kappa_{020}^{} = m_{020}^{} - m_{010}^{2},\, \kappa_{002}^{} = m_{002}^{} - m_{001}^{2}
        \end{split}
    \end{equation}
    Furthermore, substituting Eq.~(\ref{cumulant moment 1 and 2 order 3 variable}) into the related expressions obtained by comparing the coefficients, the higher-order moments can be represented using the first and second-order moments. By this means Eq.~(\ref{moment cumulant higher order 3 variable part1}) - (\ref{moment cumulant higher order 3 variable part12}) are obtained.
    \begin{equation}
    \label{moment cumulant higher order 3 variable part1}
        \begin{split}
            \displaystyle
            m_{1,1,1}^{}
            &
            = \kappa_{100}^{}\kappa_{010}^{}\kappa_{001}^{} + \kappa_{100}^{}\kappa_{011}^{} + \kappa_{010}^{}\kappa_{101}^{} + \kappa_{001}^{}\kappa_{110}^{}
            \\
            \displaystyle &
            = m_{100}^{}m_{011}^{} + m_{010}^{}m_{101}^{} + m_{001}^{}m_{110}^{}
            \\
            \displaystyle &
            -\, 2m_{100}^{}m_{010}^{}m_{001}^{}
        \end{split}
    \end{equation}
    \begin{equation}
    \label{moment cumulant higher order 3 variable part2}
        \begin{split}
            \displaystyle
            m_{1,1,2}^{}
            &
            = \kappa_{010}^{} \kappa_{100}^{} \kappa_{001}^2 + \kappa_{110}^{} \kappa_{001}^2 + 2 \kappa_{011}^{} \kappa_{100}^{} \kappa_{001}^{} 
            \\
            \displaystyle &
            +\, 2 \kappa_{010}^{} \kappa_{101}^{} \kappa_{001}^{} + \kappa_{002}^{} \kappa_{010}^{} \kappa_{100}^{} + 2 \kappa_{011}^{} \kappa_{101}^{}
            \\
            \displaystyle &
            +\, \kappa_{002}^{} \kappa_{110}^{}
            \\
            \displaystyle &
            = -2 m_{010}^{} m_{100}^{} m_{001}^2 + 2 m_{011}^{} m_{101}^{} + m_{002}^{} m_{110}^{}
        \end{split}
    \end{equation}
    \begin{equation}
    \label{moment cumulant higher order 3 variable part3}
        \begin{split}
            \displaystyle
            m_{1,1,3}^{}
            &
            = \kappa_{010}^{} \kappa_{100}^{} \kappa_{001}^3 + \kappa_{110}^{} \kappa_{001}^3 + 3 \kappa_{011}^{} \kappa_{100}^{} \kappa_{001}^2
            \\
            \displaystyle &
            +\, 3 \kappa_{010}^{} \kappa_{101}^{} \kappa_{001}^2 + 3 \kappa_{002}^{} \kappa_{010}^{} \kappa_{100}^{} \kappa_{001}^{} 
            \\
            \displaystyle &
            +\, 6 \kappa_{011}^{}\kappa_{101}^{} \kappa_{001}^{} + 3 \kappa_{002}^{} \kappa_{110}^{} \kappa_{001}^{}
            \\
            \displaystyle &
            +\, 3 \kappa_{002}^{} \kappa_{011}^{} \kappa_{100}^{} + 3 \kappa_{002}^{} \kappa_{010}^{} \kappa_{101}^{}
            \\
            \displaystyle &
            = 6 m_{010}^{} m_{100}^{} m_{001}^3 - 2 m_{110}^{} m_{001}^3
            \\
            \displaystyle &
            -\, 6 m_{011}^{} m_{100}^{} m_{001}^2 - 6 m_{010}^{} m_{101}^{} m_{001}^2
            \\
            \displaystyle &
            -\, 6 m_{002}^{} m_{010}^{} m_{100}^{} m_{001}^{} + 6m_{011}^{} m_{101}^{} m_{001}^{}
            \\
            \displaystyle &
            +\, 3 m_{002}^{} m_{110}^{} m_{001}^{} + 3 m_{002}^{} m_{011}^{} m_{100}^{}
            \\
            \displaystyle &
            +\, 3 m_{002}^{} m_{010}^{} m_{101}^{}
        \end{split}
    \end{equation}   
    \begin{equation}
    \label{moment cumulant higher order 3 variable part4}
        \begin{split}
            \displaystyle
            m_{1,1,4}^{}
            & = \kappa_{010}^{} \kappa_{100}^{} \kappa_{001}^4 + \kappa_{110}^{} \kappa_{001}^4 + 4 \kappa_{011}^{} \kappa_{100}^{} \kappa_{001}^3 
            \\
            \displaystyle &
            +\, 4 \kappa_{010}^{} \kappa_{101}^{} \kappa_{001}^3 + 6 \kappa_{002}^{} \kappa_{010}^{} \kappa_{100}^{} \kappa_{001}^2 
            \\
            \displaystyle &
            +\, 12 \kappa_{011}^{}\kappa_{101}^{} \kappa_{001}^2 + 6 \kappa_{002}^{} \kappa_{110}^{} \kappa_{001}^2 
            \\
            \displaystyle &
            +\, 12 \kappa_{002}^{} \kappa_{011}^{} \kappa_{100}^{} \kappa_{001}^{} + 12 \kappa_{002}^{} \kappa_{010}^{} \kappa_{101}^{} \kappa_{001}^{} 
            \\
            \displaystyle &
            +\, 3 \kappa_{002}^2 \kappa_{010}^{}\kappa_{100}^{} + 12 \kappa_{002}^{} \kappa_{011}^{} \kappa_{101}^{} + 3 \kappa_{002}^2 \kappa_{110}^{}
            \\
            \displaystyle &
            = 16 m_{010}^{} m_{100}^{} m_{001}^4 - 2 m_{110}^{} m_{001}^4
            \\
            \displaystyle &
            -\, 8 m_{011}^{} m_{100}^{} m_{001}^3 - 8 m_{010}^{} m_{101}^{} m_{001}^3 
            \\
            \displaystyle &
            -\, 12 m_{002}^{} m_{010}^{} m_{100}^{} m_{001}^2 + 12m_{002}^{} m_{011}^{} m_{101}^{} 
            \\
            \displaystyle &
            +\, 3 m_{002}^2 m_{110}^{}
        \end{split}
    \end{equation}   
    \begin{equation}
    \label{moment cumulant higher order 3 variable part5}
        \begin{split}
            \displaystyle
            m_{1,1,5}^{}
            &
            = \kappa_{010}^{} \kappa_{100}^{} \kappa_{001}^5 + \kappa_{110}^{} \kappa_{001}^5 + 5 \kappa_{011}^{} \kappa_{100}^{} \kappa_{001}^4 
            \\
            \displaystyle &
            +\, 5 \kappa_{010}^{} \kappa_{101}^{} \kappa_{001}^4 + 10 \kappa_{002}^{} \kappa_{010}^{} \kappa_{100}^{} \kappa_{001}^3
            \\
            \displaystyle &
            +\, 20\kappa_{011}^{} \kappa_{101}^{} \kappa_{001}^3 + 10 \kappa_{002}^{} \kappa_{110}^{} \kappa_{001}^3 
            \\
            \displaystyle &
            +\, 30 \kappa_{002}^{} \kappa_{011}^{} \kappa_{100}^{} \kappa_{001}^2 + 30 \kappa_{002}^{} \kappa_{010}^{} \kappa_{101}^{} \kappa_{001}^2
            \\
            \displaystyle &
            +\, 15 \kappa_{002}^2\kappa_{010}^{} \kappa_{100}^{} \kappa_{001}^{} + 60 \kappa_{002}^{} \kappa_{011}^{} \kappa_{1,0,1}^{} \kappa_{0,0,1}^{} 
            \\
            \displaystyle &
            +\, 15 \kappa_{002}^2 \kappa_{110}^{} \kappa_{001}^{} + 15 \kappa_{002}^2 \kappa_{011}^{} \kappa_{100}^{}
            \\
            \displaystyle &
            + 15 \kappa_{002}^2 \kappa_{010}^{}\kappa_{101}^{}
            \\
            \displaystyle &
            = -20 m_{010}^{} m_{100}^{} m_{001}^5 + 6 m_{110}^{} m_{001}^5
            \\
            \displaystyle &
            +\, 30 m_{011}^{} m_{100}^{} m_{001}^4 + 30 m_{010}^{} m_{101}^{} m_{001}^4 
            \\
            \displaystyle &
            +\, 60 m_{002}^{} m_{010}^{} m_{100}^{} m_{001}^3 - 40m_{011}^{} m_{101}^{} m_{001}^3
            \\
            \displaystyle &
            -\, 20 m_{002}^{} m_{110}^{} m_{001}^3 - 60 m_{002}^{} m_{011}^{} m_{100}^{} m_{001}^2
            \\
            \displaystyle &
            -\, 60 m_{002}^{} m_{010}^{} m_{101}^{} m_{001}^2 - 30 m_{002}^2m_{010}^{} m_{100}^{} m_{001}^{}
            \\
            \displaystyle &
            +\, 60 m_{002}^{} m_{011}^{} m_{101}^{} m_{001}^{} + 15 m_{002}^2 m_{110}^{} m_{001}^{}
            \\
            \displaystyle &
            +\, 15 m_{002}^2 m_{011}^{} m_{100}^{} + 15 m_{002}^2 m_{010}^{}m_{101}^{}
        \end{split}
    \end{equation} 
    \begin{equation}
    \label{moment cumulant higher order 3 variable part6}
        \begin{split}
            \displaystyle
            m_{1,1,6}^{}
            &
            = \kappa_{010}^{} \kappa_{100}^{} \kappa_{001}^6 + \kappa_{110}^{} \kappa_{001}^6 + 6 \kappa_{011}^{} \kappa_{100}^{} \kappa_{001}^5
            \\
            \displaystyle &
            +\, 6 \kappa_{010}^{} \kappa_{101}^{} \kappa_{001}^5 + 15 \kappa_{002}^{} \kappa_{010}^{} \kappa_{100}^{} \kappa_{001}^4 
            \\
            \displaystyle &
            +\, 30\kappa_{011}^{} \kappa_{101}^{} \kappa_{001}^4 + 15 \kappa_{002}^{} \kappa_{110}^{} \kappa_{001}^4
            \\
            \displaystyle &
            +\, 60\kappa_{002}^{} \kappa_{011}^{} \kappa_{100}^{} \kappa_{001}^3 + 60 \kappa_{002}^{} \kappa_{010}^{} \kappa_{101}^{} \kappa_{001}^3 
            \\
            \displaystyle &
            +\, 45 \kappa_{002}^2\kappa_{010}^{} \kappa_{100}^{} \kappa_{001}^2 + 180 \kappa_{002}^{} \kappa_{011}^{} \kappa_{101}^{} \kappa_{001}^2
            \\
            \displaystyle &
            +\, 45 \kappa_{002}^2 \kappa_{110}^{} \kappa_{001}^2 + 90 \kappa_{002}^2 \kappa_{011}^{} \kappa_{100}^{} \kappa_{001}^{} 
            \\
            \displaystyle &
            +\, 90\kappa_{002}^2\kappa_{010}^{} \kappa_{101}^{} \kappa_{001}^{} + 15 \kappa_{002}^3 \kappa_{010}^{} \kappa_{100}^{}
            \\
            \displaystyle &
            +\, 90 \kappa_{002}^2 \kappa_{011}^{} \kappa_{101}^{} + 15 \kappa_{002}^3 \kappa_{110}^{}
            \\
            \displaystyle &
            = -132 m_{010}^{} m_{100}^{} m_{001}^6 + 16 m_{110}^{} m_{001}^6 
            \\
            \displaystyle &
            +\, 96 m_{011}^{} m_{100}^{} m_{001}^5 + 96 m_{010}^{} m_{101}^{} m_{001}^5 
            \\
            \displaystyle &
            +\, 240 m_{002}^{} m_{010}^{} m_{100}^{}m_{001}^4 - 60 m_{011}^{} m_{101}^{} m_{001}^4 
            \\
            \displaystyle &
            -\, 30 m_{002}^{} m_{110}^{} m_{001}^4 - 120 m_{002}^{} m_{011}^{} m_{100}^{} m_{001}^3
            \\
            \displaystyle &
            -\, 120 m_{002}^{} m_{010}^{} m_{101}^{}m_{001}^3 - 90 m_{002}^2 m_{010}^{} m_{100}^{} m_{001}^2
            \\
            \displaystyle &
            +\, 90 m_{002}^2 m_{011}^{} m_{101}^{} + 15 m_{002}^3 m_{110}^{}
        \end{split}
    \end{equation}      
    \begin{equation}
    \label{moment cumulant higher order 3 variable part7}
        \begin{split}
            \displaystyle
            m_{1,1,7}^{}
            &
            = \kappa_{010}^{} \kappa_{100}^{} \kappa_{001}^7 + \kappa_{110}^{} \kappa_{001}^7 + 7 \kappa_{011}^{} \kappa_{100}^{} \kappa_{001}^6
            \\
            \displaystyle &
            +\, 7 \kappa_{010}^{} \kappa_{101}^{} \kappa_{001}^6 + 21 \kappa_{002}^{} \kappa_{010}^{} \kappa_{100}^{} \kappa_{001}^5
            \\
            \displaystyle &
            +\, 42\kappa_{011}^{} \kappa_{101}^{} \kappa_{001}^5 + 21 \kappa_{002}^{} \kappa_{110}^{} \kappa_{001}^5
            \\
            \displaystyle &
            +\, 105 \kappa_{002}^{} \kappa_{011}^{} \kappa_{100}^{} \kappa_{001}^4 + 105 \kappa_{002}^{} \kappa_{010}^{} \kappa_{101}^{} \kappa_{001}^4 
            \\
            \displaystyle &
            +\, 105 \kappa_{002}^2\kappa_{010}^{} \kappa_{100}^{} \kappa_{001}^3 + 420 \kappa_{002}^{} \kappa_{011}^{} \kappa_{101}^{} \kappa_{001}^3
            \\
            \displaystyle &
            +\, 105 \kappa_{002}^2 \kappa_{110}^{} \kappa_{001}^3 + 315 \kappa_{002}^2 \kappa_{011}^{} \kappa_{100}^{} \kappa_{001}^2 
            \\
            \displaystyle &
            +\, 315 \kappa_{002}^2 \kappa_{010}^{} \kappa_{101}^{} \kappa_{001}^2 + 105 \kappa_{002}^3 \kappa_{010}^{} \kappa_{100}^{} \kappa_{001}^{}
            \\
            \displaystyle &
            +\, 630 \kappa_{002}^2 \kappa_{011}^{} \kappa_{101}^{} \kappa_{001}^{} + 105 \kappa_{002}^3 \kappa_{110}^{} \kappa_{001}^{}
            \\
            \displaystyle &
            +\, 105\kappa_{002}^3 \kappa_{011}^{} \kappa_{100}^{} + 105 \kappa_{002}^3 \kappa_{010}^{} \kappa_{101}^{}
            \\
            \displaystyle &
            = 28 m_{010}^{} m_{100}^{} m_{001}^7 - 20 m_{110}^{} m_{001}^7 
            \\
            \displaystyle &
            -\, 140 m_{011}^{} m_{100}^{} m_{001}^6 + 630 m_{002}^2 m_{010}^{} m_{100}^{} m_{001}^3
            \\
            \displaystyle &
            -\, 420 m_{002}^{} m_{010}^{} m_{100}^{}m_{001}^5 + 252 m_{011}^{} m_{101}^{} m_{001}^5
            \\
            \displaystyle &
            +\, 126 m_{002}^{} m_{110}^{} m_{001}^5 + 630 m_{002}^{} m_{011}^{} m_{100}^{} m_{001}^4 
            \\
            \displaystyle &
            +\, 630 m_{002}^{} m_{010}^{} m_{101}^{} m_{001}^4 - 140 m_{010}^{} m_{101}^{} m_{001}^6 
            \\
            \displaystyle &
            -\, 840 m_{002}^{} m_{011}^{} m_{101}^{} m_{001}^3 - 210 m_{002}^2 m_{110}^{} m_{001}^3
            \\
            \displaystyle &
            -\, 630 m_{002}^2 m_{011}^{}m_{100}^{} m_{001}^2 + 105 m_{002}^3 m_{011}^{} m_{100}^{} 
            \\
            \displaystyle &
            -\, 210 m_{002}^3 m_{010}^{} m_{100}^{} m_{001}^{} + 105 m_{002}^3 m_{010}^{} m_{101}^{} 
            \\
            \displaystyle &
            +\, 105 m_{002}^3 m_{110}^{} m_{001}^{} - 630 m_{002}^2 m_{010}^{} m_{101}^{} m_{001}^2 
            \\
            \displaystyle &
            +\, 630 m_{002}^2 m_{011}^{} m_{101}^{} m_{001}^{}
        \end{split}
    \end{equation}
     \begin{equation}
    \label{moment cumulant higher order 3 variable part8}
        \begin{split}
            \displaystyle
            m_{1,1,8}^{}
            &
            = \kappa_{010}^{} \kappa_{100}^{} \kappa_{001}^8 + \kappa_{110}^{} \kappa_{001}^8 + 8 \kappa_{011}^{} \kappa_{100}^{} \kappa_{001}^7 
            \\
            \displaystyle &
            +\, 8 \kappa_{010}^{} \kappa_{101}^{} \kappa_{001}^7 + 28 \kappa_{002}^{} \kappa_{010}^{} \kappa_{100}^{} \kappa_{001}^6 
            \\
            \displaystyle &
            +\, 56 \kappa_{011}^{} \kappa_{101}^{} \kappa_{001}^6 + 28 \kappa_{002}^{} \kappa_{110}^{} \kappa_{001}^6
            \\
            \displaystyle &
            +\, 168 \kappa_{002}^{} \kappa_{011}^{} \kappa_{100}^{} \kappa_{001}^5 + 168 \kappa_{002}^{} \kappa_{010}^{} \kappa_{101}^{} \kappa_{001}^5 
            \\
            \displaystyle &
            +\, 210 \kappa_{002}^2 \kappa_{010}^{} \kappa_{100}^{} \kappa_{001}^4 + 840 \kappa_{002}^{} \kappa_{011}^{} \kappa_{101}^{} \kappa_{001}^4
            \\
            \displaystyle &
            +\, 210 \kappa_{002}^2 \kappa_{110}^{} \kappa_{001}^4 + 840 \kappa_{002}^2 \kappa_{011}^{} \kappa_{100}^{} \kappa_{001}^3 
            \\
            \displaystyle &
            +\, 840 \kappa_{002}^2 \kappa_{010}^{} \kappa_{101}^{} \kappa_{001}^3 + 420 \kappa_{002}^3 \kappa_{010}^{} \kappa_{100}^{} \kappa_{001}^2
            \\
            \displaystyle &
            +\, 2520 \kappa_{002}^2 \kappa_{011}^{} \kappa_{101}^{} \kappa_{001}^2 + 420 \kappa_{002}^3 \kappa_{110}^{} \kappa_{001}^2 
            \\
            \displaystyle &
            +\, 840 \kappa_{002}^3 \kappa_{011}^{} \kappa_{100}^{} \kappa_{001}^{} + 840 \kappa_{002}^3 \kappa_{010}^{} \kappa_{101}^{} \kappa_{001}^{}
            \\
            \displaystyle &
            +\, 105 \kappa_{002}^4 \kappa_{010}^{} \kappa_{100}^{} + 840 \kappa_{002}^3 \kappa_{011}^{} \kappa_{101}^{}
            \\
            \displaystyle &
            +\, 105 \kappa_{002}^4 \kappa_{110}^{}
            \\
            \displaystyle &
            = 1216 m_{010}^{} m_{100}^{} m_{001}^8 - 1056 m_{011}^{} m_{100}^{} m_{001}^7
            \\
            \displaystyle &
            -\, 132 m_{110}^{} m_{001}^8 + 3360 m_{002}^2 m_{010}^{} m_{100}^{} m_{001}^4 
            \\
            \displaystyle &
            -\, 3696 m_{002}^{} m_{010}^{} m_{100}^{} m_{001}^6 + 896 m_{011}^{} m_{101}^{} m_{001}^6
            \\
            \displaystyle &
            +\, 448 m_{002}^{} m_{110}^{} m_{001}^6 + 2688 m_{002}^{} m_{011}^{} m_{100}^{} m_{001}^5
            \\
            \displaystyle &
            +\, 2688 m_{002}^{} m_{010}^{} m_{101}^{} m_{001}^5 - 1056 m_{010}^{} m_{101}^{} m_{001}^7 
            \\
            \displaystyle &
            -\, 1680 m_{002}^{} m_{011}^{} m_{101}^{} m_{001}^4 - 420 m_{002}^2 m_{110}^{} m_{001}^4
            \\
            \displaystyle &
            -\, 1680 m_{002}^2 m_{011}^{} m_{100}^{} m_{001}^3 + 105 m_{002}^4 m_{110}^{}
            \\
            \displaystyle &
            -\, 840 m_{002}^3 m_{010}^{} m_{100}^{} m_{001}^2 + 840 m_{002}^3 m_{011}^{} m_{101}^{} 
            \\
            \displaystyle &
            -\, 1680 m_{002}^2 m_{010}^{} m_{101}^{} m_{001}^3
        \end{split}
    \end{equation}   
    \begin{equation}
    \label{moment cumulant higher order 3 variable part9}
        \begin{split}
            \displaystyle &
            m_{1,1,9}^{}
            \\
            \displaystyle &
            = \kappa_{010}^{} \kappa_{100}^{} \kappa_{001}^9 + \kappa_{110}^{} \kappa_{001}^9 + 9 \kappa_{011}^{} \kappa_{100}^{} \kappa_{001}^8
            \\
            \displaystyle &
            +\, 9 \kappa_{010}^{} \kappa_{101}^{}\kappa_{001}^8 + 36 \kappa_{002}^{} \kappa_{010}^{} \kappa_{100}^{} \kappa_{001}^7 
            \\
            \displaystyle &
            +\, 72 \kappa_{011}^{} \kappa_{101}^{} \kappa_{001}^7 + 252 \kappa_{002}^{} \kappa_{011}^{} \kappa_{100}^{} \kappa_{001}^6
            \\
            \displaystyle &
            +\, 36 \kappa_{002}^{} \kappa_{110}^{} \kappa_{001}^7 + 252 \kappa_{002}^{} \kappa_{010}^{} \kappa_{101}^{} \kappa_{001}^6 
            \\
            \displaystyle &
            +\, 378 \kappa_{002}^2 \kappa_{010}^{} \kappa_{100}^{} \kappa_{001}^5 + 1512 \kappa_{002}^{} \kappa_{011}^{} \kappa_{101}^{} \kappa_{001}^5
            \\
            \displaystyle &
            +\, 378 \kappa_{002}^2 \kappa_{110}^{} \kappa_{001}^5 + 1890 \kappa_{002}^2 \kappa_{011}^{} \kappa_{100}^{} \kappa_{001}^4
            \\
            \displaystyle &
            +\, 1890 \kappa_{002}^2 \kappa_{010}^{} \kappa_{101}^{} \kappa_{001}^4 + 1260 \kappa_{002}^3 \kappa_{010}^{} \kappa_{100}^{} \kappa_{001}^3 
            \\
            \displaystyle &
            +\, 7560 \kappa_{002}^2 \kappa_{011}^{} \kappa_{101}^{} \kappa_{001}^3 + 1260 \kappa_{002}^3 \kappa_{110}^{} \kappa_{001}^3
            \\
            \displaystyle &
            +\, 3780 \kappa_{002}^3 \kappa_{011}^{} \kappa_{100}^{} \kappa_{001}^2 + 3780 \kappa_{002}^3 \kappa_{010}^{} \kappa_{101}^{} \kappa_{001}^2 
            \\
            \displaystyle &
            +\, 945 \kappa_{002}^4 \kappa_{010}^{} \kappa_{100}^{} \kappa_{001}^{} + 7560 \kappa_{002}^3 \kappa_{011}^{} \kappa_{101}^{} \kappa_{001}^{} 
            \\
            \displaystyle &
            +\, 945 \kappa_{002}^4 \kappa_{110}^{} \kappa_{001}^{} + 945 \kappa_{002}^4 \kappa_{011}^{} \kappa_{100}^{}
            \\
            \displaystyle &
            +\, 945 \kappa_{002}^4 \kappa_{010}^{} \kappa_{101}^{}
            \\
            \displaystyle &
            = 936 m_{010}^{} m_{100}^{} m_{001}^9 + 7560 m_{002}^3 m_{010}^{} m_{100}^{} m_{001}^3
            \\
            \displaystyle &
            +\, 252 m_{010}^{} m_{101}^{} m_{001}^8 + 1008 m_{002}^{} m_{010}^{} m_{100}^{} m_{001}^7 
            \\
            \displaystyle &
            -\, 1440 m_{011}^{} m_{101}^{} m_{001}^7 - 5040 m_{002}^{} m_{010}^{} m_{101}^{} m_{001}^6
            \\
            \displaystyle &
            -\, 5040 m_{002}^{} m_{011}^{} m_{100}^{} m_{001}^6 - 720 m_{002}^{} m_{110}^{} m_{001}^7 
            \\
            \displaystyle &
            -\, 7560 m_{002}^2 m_{010}^{} m_{100}^{} m_{001}^5 + 945 m_{002}^4 m_{011}^{} m_{100}^{}
            \\
            \displaystyle &
            +\, 2268 m_{002}^2 m_{110}^{} m_{001}^5 + 11340 m_{002}^2 m_{011}^{} m_{100}^{} m_{001}^4 
            \\
            \displaystyle &
            +\, 11340 m_{002}^2 m_{010}^{} m_{101}^{} m_{001}^4 + 252 m_{011}^{} m_{100}^{} m_{001}^8
            \\
            \displaystyle &
            -\, 15120 m_{002}^2 m_{011}^{} m_{101}^{} m_{001}^3 - 2520 m_{002}^3 m_{110}^{} m_{001}^3
            \\
            \displaystyle &
            -\, 7560 m_{002}^3 m_{011}^{} m_{100}^{} m_{001}^2 + 945 m_{002}^4 m_{010}^{} m_{101}^{}
            \\
            \displaystyle &
            -\, 1890 m_{002}^4 m_{010}^{} m_{100}^{} m_{001}^{} + 28 m_{110}^{} m_{001}^9
            \\
            \displaystyle &
            +\, 945 m_{002}^4 m_{110}^{} m_{001}^{} + 9072 m_{002}^{} m_{011}^{} m_{101}^{} m_{001}^5
            \\
            \displaystyle &
            +\, 7560 m_{002}^3 m_{011}^{} m_{101}^{} m_{001}^{} - 7560 m_{002}^3 m_{010}^{} m_{101}^{} m_{001}^2 
        \end{split}
    \end{equation}   
    \begin{equation}
    \label{moment cumulant higher order 3 variable part10}
        \begin{split}
            \displaystyle &
            m_{1,1,10}^{}
            \\
            \displaystyle &
            = \kappa_{010}^{} \kappa_{100}^{} \kappa_{001}^{10} + \kappa_{110}^{} \kappa_{001}^{10} + 10 \kappa_{011}^{} \kappa_{100}^{} \kappa_{001}^9
            \\
            \displaystyle &
            +\, 10 \kappa_{010}^{} \kappa_{101}^{} \kappa_{001}^9 + 90 \kappa_{011}^{} \kappa_{101}^{} \kappa_{001}^8 
            \\
            \displaystyle &
            +\, 45 \kappa_{002}^{} \kappa_{010}^{} \kappa_{100}^{} \kappa_{001}^8  + 360 \kappa_{002}^{} \kappa_{011}^{} \kappa_{100}^{} \kappa_{001}^7
            \\
            \displaystyle &
            +\, 360 \kappa_{002}^{} \kappa_{010}^{} \kappa_{101}^{} \kappa_{001}^7 + 630 \kappa_{002}^2 \kappa_{010}^{} \kappa_{100}^{} \kappa_{001}^6
            \\
            \displaystyle &
            +\, 2520 \kappa_{002}^{} \kappa_{011}^{} \kappa_{101}^{} \kappa_{001}^6 + 630 \kappa_{002}^2 \kappa_{110}^{} \kappa_{001}^6 
            \\
            \displaystyle &
            +\, 3780 \kappa_{002}^2 \kappa_{011}^{} \kappa_{100}^{} \kappa_{001}^5 + 3780 \kappa_{002}^2 \kappa_{010}^{} \kappa_{101}^{} \kappa_{001}^5
            \\
            \displaystyle &
            +\, 3150 \kappa_{002}^3 \kappa_{010}^{} \kappa_{100}^{} \kappa_{001}^4 + 18900 \kappa_{002}^2 \kappa_{011}^{} \kappa_{101}^{} \kappa_{001}^4
            \\
            \displaystyle &
            +\, 3150 \kappa_{002}^3 \kappa_{110}^{} \kappa_{001}^4 + 12600 \kappa_{002}^3 \kappa_{011}^{} \kappa_{100}^{} \kappa_{001}^3 
            \\
            \displaystyle &
            +\, 12600 \kappa_{002}^3 \kappa_{010}^{} \kappa_{101}^{} \kappa_{001}^3 + 4725 \kappa_{002}^4 \kappa_{010}^{} \kappa_{100}^{} \kappa_{001}^2
            \\
            \displaystyle &
            +\, 37800 \kappa_{002}^3 \kappa_{011}^{} \kappa_{101}^{} \kappa_{001}^2 + 4725 \kappa_{002}^4 \kappa_{110}^{} \kappa_{001}^2
            \\
            \displaystyle &
            +\, 9450 \kappa_{002}^4 \kappa_{011}^{} \kappa_{100}^{} \kappa_{001}^{} + 9450 \kappa_{002}^4 \kappa_{010}^{} \kappa_{101}^{} \kappa_{001}^{} 
            \\
            \displaystyle &
            +\, 945 \kappa_{002}^5 \kappa_{010}^{} \kappa_{100}^{} + 9450 \kappa_{002}^4 \kappa_{011}^{} \kappa_{101}^{}
            \\
            \displaystyle &
            +\, 945 \kappa_{002}^5 \kappa_{110}^{} + 45 \kappa_{002}^{} \kappa_{110}^{} \kappa_{001}^8
            \\
            \displaystyle &
            = -12440 m_{010}^{} m_{100}^{} m_{001}^{10} + 1216 m_{110}^{} m_{001}^{10}
            \\
            \displaystyle &
            +\, 12160 m_{011}^{} m_{100}^{} m_{001}^9 - 37800 m_{002}^2 m_{011}^{} m_{101}^{} m_{001}^4 
            \\
            \displaystyle &
            +\, 54720 m_{002}^{} m_{010}^{} m_{100}^{} m_{001}^8 - 11880 m_{011}^{} m_{101}^{} m_{001}^8
            \\
            \displaystyle &
            -\, 5940 m_{002}^{} m_{110}^{} m_{001}^8 - 47520 m_{002}^{} m_{011}^{} m_{100}^{} m_{001}^7
            \\
            \displaystyle &
            -\, 47520 m_{002}^{} m_{010}^{} m_{101}^{} m_{001}^7 + 9450 m_{002}^4 m_{011}^{} m_{101}^{}
            \\
            \displaystyle &
            +\, 40320 m_{002}^{} m_{011}^{} m_{101}^{} m_{001}^6 + 10080 m_{002}^2 m_{110}^{} m_{001}^6 
            \\
            \displaystyle &
            +\, 60480 m_{002}^2 m_{011}^{} m_{100}^{} m_{001}^5 + 945 m_{002}^5 m_{110}^{}
            \\
            \displaystyle &
            +\, 50400 m_{002}^3 m_{010}^{} m_{100}^{} m_{001}^4 + 12160 m_{010}^{} m_{101}^{} m_{001}^9 
            \\
            \displaystyle &
            -\, 6300 m_{002}^3 m_{110}^{} m_{001}^4 - 25200 m_{002}^3 m_{011}^{} m_{100}^{} m_{001}^3 
            \\
            \displaystyle &
            -\, 25200 m_{002}^3 m_{010}^{} m_{101}^{} m_{001}^3 - 9450 m_{002}^4 m_{010}^{} m_{100}^{} m_{001}^2
            \\
            \displaystyle &
            -\, 83160 m_{002}^2 m_{010}^{} m_{100}^{} m_{001}^6 + 60480 m_{002}^2 m_{010}^{} m_{101}^{} m_{001}^5
        \end{split}
    \end{equation}   
    \begin{equation}
    \label{moment cumulant higher order 3 variable part11}
        \begin{split}
            \displaystyle &
            m_{1,1,11}^{}
            \\
            \displaystyle &
            = \kappa_{010}^{} \kappa_{100}^{} \kappa_{001}^{11} + \kappa_{110}^{} \kappa_{001}^{11} + 11 \kappa_{011}^{} \kappa_{100}^{} \kappa_{001}^{10}
            \\
            \displaystyle &
            +\, 11 \kappa_{010}^{} \kappa_{101}^{} \kappa_{001}^{10} + 55 \kappa_{002}^{} \kappa_{010}^{} \kappa_{100}^{} \kappa_{001}^9 
            \\
            \displaystyle &
            +\, 110 \kappa_{011}^{} \kappa_{101}^{} \kappa_{001}^9 + 495 \kappa_{002}^{} \kappa_{011}^{} \kappa_{100}^{} \kappa_{001}^8 
            \\
            \displaystyle &
            +\, 55 \kappa_{002}^{} \kappa_{110}^{} \kappa_{001}^9 + 495 \kappa_{002}^{} \kappa_{010}^{} \kappa_{101}^{} \kappa_{001}^8 
            \\
            \displaystyle &
            +\, 990 \kappa_{002}^2 \kappa_{010}^{} \kappa_{100}^{} \kappa_{001}^7 + 3960 \kappa_{002}^{} \kappa_{011}^{} \kappa_{101}^{} \kappa_{001}^7
            \\
            \displaystyle &
            +\, 990 \kappa_{002}^2 \kappa_{110}^{} \kappa_{001}^7 + 6930 \kappa_{002}^2 \kappa_{011}^{} \kappa_{100}^{} \kappa_{001}^6
            \\
            \displaystyle &
            +\, 6930 \kappa_{002}^2 \kappa_{010}^{} \kappa_{101}^{} \kappa_{001}^6 + 6930 \kappa_{002}^3 \kappa_{010}^{} \kappa_{100}^{} \kappa_{001}^5 
            \\
            \displaystyle &
            +\, 41580 \kappa_{002}^2 \kappa_{011}^{} \kappa_{101}^{} \kappa_{001}^5 + 6930 \kappa_{002}^3 \kappa_{110}^{} \kappa_{001}^5 
            \\
            \displaystyle &
            +\, 34650 \kappa_{002}^3 \kappa_{011}^{} \kappa_{100}^{} \kappa_{001}^4 + 34650 \kappa_{002}^3 \kappa_{010}^{} \kappa_{101}^{} \kappa_{001}^4
            \\
            \displaystyle &
            +\, 17325 \kappa_{002}^4 \kappa_{010}^{} \kappa_{100}^{} \kappa_{001}^3 + 138600 \kappa_{002}^3 \kappa_{011}^{} \kappa_{101}^{} \kappa_{001}^3 
            \\
            \displaystyle &
            +\, 17325 \kappa_{002}^4 \kappa_{110}^{} \kappa_{001}^3 + 51975 \kappa_{002}^4 \kappa_{011}^{} \kappa_{100}^{} \kappa_{001}^2
            \\
            \displaystyle &
            +\, 51975 \kappa_{002}^4 \kappa_{010}^{} \kappa_{101}^{} \kappa_{001}^2 + 10395 \kappa_{002}^5 \kappa_{010}^{} \kappa_{100}^{} \kappa_{001}^{} 
            \\
            \displaystyle &
            +\, 103950 \kappa_{002}^4 \kappa_{011}^{} \kappa_{101}^{} \kappa_{001}^{} + 10395 \kappa_{002}^5 \kappa_{110}^{} \kappa_{001}^{} 
            \\
            \displaystyle &
            +\, 10395 \kappa_{002}^5 \kappa_{011}^{} \kappa_{100}^{} + 10395 \kappa_{002}^5 \kappa_{010}^{} \kappa_{101}^{}
            \\
            \displaystyle &
            = -23672 m_{010}^{} m_{100}^{} m_{001}^{11} - 20790 m_{002}^5 m_{010}^{} m_{100}^{} m_{001}^{}
            \\
            \displaystyle &
            +\, 936 m_{110}^{} m_{001}^{11} + 27720 m_{002}^2 m_{010}^{} m_{100}^{} m_{001}^7
            \\
            \displaystyle &
            +\, 51480 m_{002}^{} m_{010}^{} m_{100}^{} m_{001}^9 + 3080 m_{011}^{} m_{101}^{} m_{001}^9
            \\
            \displaystyle &
            +\, 1540 m_{002}^{} m_{110}^{} m_{001}^9 + 13860 m_{002}^{} m_{011}^{} m_{100}^{} m_{001}^8
            \\
            \displaystyle &
            +\, 13860 m_{002}^{} m_{010}^{} m_{101}^{} m_{001}^8 + 10296 m_{010}^{} m_{101}^{} m_{001}^{10} 
            \\
            \displaystyle &
            -\, 79200 m_{002}^{} m_{011}^{} m_{101}^{} m_{001}^7 - 19800 m_{002}^2 m_{110}^{} m_{001}^7 
            \\
            \displaystyle &
            -\, 138600 m_{002}^2 m_{011}^{} m_{100}^{} m_{001}^6 + 10395 m_{002}^5 m_{010}^{} m_{101}^{}
            \\
            \displaystyle &
            -\, 138600 m_{002}^3 m_{010}^{} m_{100}^{} m_{001}^5 + 10395 m_{002}^5 m_{011}^{} m_{100}^{}
            \\
            \displaystyle &
            +\, 41580 m_{002}^3 m_{110}^{} m_{001}^5 + 207900 m_{002}^3 m_{011}^{} m_{100}^{} m_{001}^4
            \\
            \displaystyle &
            -\, 277200 m_{002}^3 m_{011}^{} m_{101}^{} m_{001}^3 - 34650 m_{002}^4 m_{110}^{} m_{001}^3
            \\
            \displaystyle &
            +\, 10395 m_{002}^5 m_{110}^{} m_{001}^{} + 249480 m_{002}^2 m_{011}^{} m_{101}^{} m_{001}^5 
            \\
            \displaystyle &
            +\, 207900 m_{002}^3 m_{010}^{} m_{101}^{} m_{001}^4 + 10296 m_{011}^{} m_{100}^{} m_{001}^{10}
            \\
            \displaystyle &
            +\, 103950 m_{002}^4 m_{010}^{} m_{100}^{} m_{001}^3 - 138600 m_{002}^2 m_{010}^{} m_{101}^{} m_{001}^6
            \\
            \displaystyle &
            -\, 103950 m_{002}^4 m_{011}^{} m_{100}^{} m_{001}^2- 103950 m_{002}^4 m_{010}^{} m_{101}^{} m_{001}^2
            \\
            \displaystyle &
            +\, 103950 m_{002}^4 m_{011}^{} m_{101}^{} m_{001}^{}
        \end{split}
    \end{equation}
    \begin{equation}
    \label{moment cumulant higher order 3 variable part12}
        \begin{split}
            \displaystyle &
            m_{1,1,12}^{}
            \\
            \displaystyle &
            = \kappa_{010}^{} \kappa_{100}^{} \kappa_{001}^{12} + 12 \kappa_{011}^{} \kappa_{100}^{} \kappa_{001}^{11} + 66 \kappa_{002}^{} \kappa_{110}^{} \kappa_{001}^{10}
            \\
            \displaystyle &
            +\, 66 \kappa_{002}^{} \kappa_{010}^{} \kappa_{100}^{} \kappa_{001}^{10} + 132 \kappa_{011}^{} \kappa_{101}^{} \kappa_{001}^{10} + \kappa_{110}^{} \kappa_{001}^{12} 
            \\
            \displaystyle &
            +\, 660 \kappa_{002}^{} \kappa_{011}^{} \kappa_{100}^{} \kappa_{001}^9 + 660 \kappa_{002}^{} \kappa_{010}^{} \kappa_{101}^{} \kappa_{001}^9 
            \\
            \displaystyle &
            +\, 1485 \kappa_{002}^2 \kappa_{010}^{} \kappa_{100}^{} \kappa_{001}^8 + 5940 \kappa_{002}^{} \kappa_{011}^{} \kappa_{101}^{} \kappa_{001}^8 
            \\
            \displaystyle &
            +\, 1485 \kappa_{002}^2 \kappa_{110}^{} \kappa_{001}^8 + 11880 \kappa_{002}^2 \kappa_{011}^{} \kappa_{100}^{} \kappa_{001}^7 
            \\
            \displaystyle &
            +\, 11880 \kappa_{002}^2 \kappa_{010}^{} \kappa_{101}^{} \kappa_{001}^7 + 13860 \kappa_{002}^3 \kappa_{010}^{} \kappa_{100}^{} \kappa_{001}^6 
            \\
            \displaystyle &
            +\, 83160 \kappa_{002}^2 \kappa_{011}^{} \kappa_{101}^{} \kappa_{001}^6 + 13860 \kappa_{002}^3 \kappa_{110}^{} \kappa_{001}^6 
            \\
            \displaystyle &
            +\, 83160 \kappa_{002}^3 \kappa_{011}^{} \kappa_{100}^{} \kappa_{001}^5 + 83160 \kappa_{002}^3 \kappa_{010}^{} \kappa_{101}^{} \kappa_{001}^5 
            \\
            \displaystyle &
            +\, 51975 \kappa_{002}^4 \kappa_{010}^{} \kappa_{100}^{} \kappa_{001}^4 + 415800 \kappa_{002}^3 \kappa_{011}^{} \kappa_{101}^{} \kappa_{001}^4 
            \\
            \displaystyle &
            +\, 51975 \kappa_{002}^4 \kappa_{110}^{} \kappa_{001}^4 + 207900 \kappa_{002}^4 \kappa_{011}^{} \kappa_{100}^{} \kappa_{001}^3 
            \\
            \displaystyle &
            +\, 207900 \kappa_{002}^4 \kappa_{010}^{} \kappa_{101}^{} \kappa_{001}^3 + 62370 \kappa_{002}^5 \kappa_{010}^{} \kappa_{100}^{} \kappa_{001}^2
            \\
            \displaystyle &
            +\, 623700 \kappa_{002}^4 \kappa_{011}^{} \kappa_{101}^{} \kappa_{001}^2 + 62370 \kappa_{002}^5 \kappa_{110}^{} \kappa_{001}^2 
            \\
            \displaystyle &
            +\, 124740 \kappa_{002}^5 \kappa_{011}^{} \kappa_{100}^{} \kappa_{001}^{} + 124740 \kappa_{002}^5 \kappa_{010}^{} \kappa_{101}^{} \kappa_{001}^{} 
            \\
            \displaystyle &
            +\, 10395 \kappa_{002}^6 \kappa_{010}^{} \kappa_{100}^{} + 124740 \kappa_{002}^5 \kappa_{011}^{} \kappa_{101}^{}
            \\
            \displaystyle &
            +\, 10395 \kappa_{002}^6 \kappa_{110}^{} + 12 \kappa_{010}^{} \kappa_{101}^{} \kappa_{001}^{11}
            \\
            \displaystyle &
            = - 12440 m_{110}^{} m_{001}^{12} + 1805760 m_{002}^2 m_{010}^{} m_{100}^{} m_{001}^8
            \\
            \displaystyle &
            -\, 149280 m_{011}^{} m_{100}^{} m_{001}^{11} - 1568160 m_{002}^2 m_{010}^{} m_{101}^{} m_{001}^7 
            \\
            \displaystyle &
            -\, 821040 m_{002}^{} m_{010}^{} m_{100}^{} m_{001}^{10} + 160512 m_{011}^{} m_{101}^{} m_{001}^{10} 
            \\
            \displaystyle &
            +\, 80256 m_{002}^{} m_{110}^{} m_{001}^{10} + 802560 m_{002}^{} m_{011}^{} m_{100}^{} m_{001}^9
            \\
            \displaystyle &
            +\, 802560 m_{002}^{} m_{010}^{} m_{101}^{} m_{001}^9 + 138048 m_{010}^{} m_{100}^{} m_{001}^{12}
            \\
            \displaystyle &
            -\, 784080 m_{002}^{} m_{011}^{} m_{101}^{} m_{001}^8 - 196020 m_{002}^2 m_{110}^{} m_{001}^8
            \\
            \displaystyle &
            -\, 1568160 m_{002}^2 m_{011}^{} m_{100}^{} m_{001}^7 - 149280 m_{010}^{} m_{101}^{} m_{001}^{11}
            \\
            \displaystyle &
            -\, 1829520 m_{002}^3 m_{010}^{} m_{100}^{} m_{001}^6 + 10395 m_{002}^6 m_{110}^{}
            \\
            \displaystyle &
            +\, 221760 m_{002}^3 m_{110}^{} m_{001}^6 + 1330560 m_{002}^3 m_{011}^{} m_{100}^{} m_{001}^5 
            \\
            \displaystyle &
            +\, 1330560 m_{002}^3 m_{010}^{} m_{101}^{} m_{001}^5 + 831600 m_{002}^4 m_{010}^{} m_{100}^{} m_{001}^4 
            \\
            \displaystyle &
            -\, 831600 m_{002}^3 m_{011}^{} m_{101}^{} m_{001}^4 - 103950 m_{002}^4 m_{110}^{} m_{001}^4
            \\
            \displaystyle &
            -\, 415800 m_{002}^4 m_{011}^{} m_{100}^{} m_{001}^3 - 415800 m_{002}^4 m_{010}^{} m_{101}^{} m_{001}^3
            \\
            \displaystyle &
            -\, 124740 m_{002}^5 m_{010}^{} m_{100}^{} m_{001}^2 + 124740 m_{002}^5 m_{011}^{} m_{101}^{} 
            \\
            \displaystyle &
            +\, 1330560 m_{002}^2 m_{011}^{} m_{101}^{} m_{001}^6
        \end{split}
    \end{equation}